\documentclass[12pt]{amsart}
\usepackage{amssymb,amscd}
\usepackage{pstricks}

\newtheorem{theorem}{Theorem}[section]
\newtheorem{lemma}[theorem]{Lemma}
\newtheorem{proposition}[theorem]{Proposition}

\newtheorem{definition}[theorem]{Definition}
\newtheorem{conj}[theorem]{Conjecture}
\numberwithin{equation}{section}

\newcommand{\Z}{\mathbb{Z}}
\newcommand{\R}{\mathbb{R}}
\newcommand{\C}{\mathbb{C}}
\newcommand{\CP}{\mathbb{CP}}
\newcommand{\PP}{\mathbb{P}}
\renewcommand{\Re}{\mathrm{Re}\,}
\renewcommand{\Im}{\mathrm{Im}\,}
\renewcommand{\L}{\mathcal{L}}
\newcommand{\m}{\mathfrak{m}}
\newcommand{\hol}{\mathrm{hol}}
\newcommand{\F}{\mathbb{F}}

\title{Special Lagrangian fibrations, wall-crossing,
and mirror symmetry}
\author{Denis Auroux}
\address{Department of Mathematics, M.I.T., Cambridge MA 02139, USA}
\email{auroux@math.mit.edu}
\thanks{
This work was partially supported by NSF grants DMS-0600148 and 
DMS-0652630.}

\begin{document} 
\begin{abstract} 
In this survey paper, we briefly review
various aspects of the SYZ approach to
mirror symmetry for non-Calabi-Yau varieties, focusing in
particular on Lagrangian fibrations and wall-crossing phenomena in Floer 
homology. Various examples are presented, some of them new.
\end{abstract}

\maketitle
\tableofcontents

\section{Introduction}

While mirror symmetry first arose as a set of predictions relating
Hodge structures and quantum cohomology for Calabi-Yau 3-folds (see e.g.\
\cite{CDGP,CoxKatz}), it has since been extended in spectacular ways.
To mention just a few key advances, Kontsevich's homological mirror
conjecture \cite{KoICM} has recast mirror symmetry in the language of
derived categories of coherent sheaves and Fukaya categories;
the Strominger-Yau-Zaslow (SYZ) conjecture \cite{SYZ} has provided the basis
for a geometric understanding of mirror symmetry; and mirror symmetry
has been extended beyond the Calabi-Yau setting, by considering
Landau-Ginzburg models (see e.g.\ \cite{HV,KoENS}).

In this paper, we briefly discuss various aspects of mirror symmetry
from the perspective of Lagrangian torus fibrations, i.e.\ following
the Strominger-Yau-Zaslow philosophy \cite{SYZ}. We mostly focus on
the case of K\"ahler manifolds with effective anticanonical divisors, 
along the same general lines as \cite{Au}. The two main phenomena
that we would like to focus on here are, on one hand, wall-crossing
in Floer homology and its role in determining ``instanton corrections''
to the complex geometry of the mirror; and on the other hand, the
possibility of ``transferring'' mirror symmetry from a given K\"ahler 
manifold to a Calabi-Yau submanifold. 

The paper is essentially expository in
nature, expanding on the themes already present in \cite{Au}.
The discussion falls far short of the level of
sophistication present in the works of Kontsevich-Soibelman \cite{KS1,KS2},
Gross-Siebert \cite{GS1,GS2}, or Fukaya-Oh-Ohta-Ono \cite{FO3book,FO3toric};
rather, our goal is to show how various important
ideas in the modern understanding of mirror symmetry naturally arise
from the perspective of a symplectic geometer, and to illustrate them
by simple examples. Accordingly, most of the
results mentioned here are not new, though to our knowledge some of
them have not appeared anywhere in the literature.

Another word of warning is in order:
we have swept under the rug many of the issues related to the rigorous
construction of Lagrangian Floer theory, and generally speaking we take 
an optimistic view of issues such as
the existence of fundamental chains for moduli spaces of discs
and the convergence of various Floer-theoretic quantities. These happen
not to be issues in the examples we consider, but can be serious obstacles
in the general case.

The rest of this paper is organized as follows: in Section \ref{s:background}
we review the SYZ approach to the construction of mirror
pairs, and the manner in which the mirror superpotential arises naturally
as a Floer-theoretic obstruction in the non Calabi-Yau case. Section
\ref{s:examples} presents various elementary examples, focusing on
wall-crossing phenomena and instanton corrections. Section \ref{s:floer}
discusses some issues related to convergent power series Floer homology.
Finally, Section \ref{s:relms} focuses on mirror
symmetry in the {\it relative} setting, namely for a Calabi-Yau 
hypersurface representing the anticanonical class inside a K\"ahler
manifold, or more generally for a complete intersection.%
\medskip

\noindent {\bf Acknowledgements.} The ideas presented here were influenced 
in a decisive manner by numerous discussions with Mohammed Abouzaid, Paul
Seidel, and Ludmil Katzarkov. Some of the topics presented here also owe
a lot to conversations with Dima Orlov, Mark Gross, and Kenji Fukaya.
Finally, I am grateful to Mohammed Abouzaid for valuable comments on
the exposition.
This work was partially supported by NSF grants DMS-0600148 and 
DMS-0652630.

\section{Lagrangian tori and mirror symmetry} \label{s:background}

\subsection{Lagrangian tori and the SYZ conjecture} \label{ss:syz}
The SYZ conjecture essentially asserts that mirror pairs of Calabi-Yau
manifolds should carry dual special Lagrangian torus fibrations
\cite{SYZ}. This statement should be understood with
suitable qualifiers (near the large complex structure limit, with
instanton corrections, etc.), but it nonetheless gives the basic
template for the geometric construction of mirror pairs.

  From this perspective, to construct the mirror of a given Calabi-Yau
manifold $X$, one should first try to construct a special Lagrangian torus
fibration $f:X\to B$. This is a difficult problem, but assuming
it has been solved, the first guess for the mirror manifold $X^\vee$
is then the total space of the dual fibration $f^\vee$. 
Given a torus $T$, the dual
torus $T^\vee=\mathrm{Hom}(\pi_1(T),S^1)$ can be viewed as a moduli space
of rank 1 unitary local systems (i.e., flat unitary connections up to
gauge equivalence) over $T$; hence, points of the dual
fibration parametrize pairs consisting of a special Lagrangian fiber
in $X$ and a unitary local system over it.
\medskip

More precisely, let $(X,J,\omega)$ be a K\"ahler manifold of complex
dimension $n$, equipped 
with a nonvanishing holomorphic volume form $\Omega\in \Omega^{n,0}(X)$.
This is sometimes called an ``almost Calabi-Yau'' manifold (to distinguish
it from a genuine Calabi-Yau, where one would also require the 
norm of $\Omega$ with respect to the K\"ahler metric to be constant).
It is an elementary fact that the restriction of $\Omega$ to a Lagrangian
submanifold $L\subset X$ is a nowhere vanishing complex-valued $n$-form.

\begin{definition}
A Lagrangian submanifold $L\subset X$ is {\em special Lagrangian}
if the argument of $\Omega_{|L}$ is constant.
\end{definition}

The value of the constant depends only on the homology class of $L$,
and we will usually normalize $\Omega$ so that it is a multiple
of $\pi/2$. For simplicity, in the rest of this paragraph we will assume
that $\Omega_{|L}$ is a positive real multiple of the real volume form
$\mathrm{vol}_g$ induced by the K\"ahler metric $g=\omega(\cdot,J\cdot)$.

The following classical result is due to McLean \cite{mclean} (in the
Calabi-Yau setting; see \S 9 of \cite{joycenotes} or Proposition 2.5 
of \cite{Au} for the almost Calabi-Yau case):

\begin{proposition}[McLean]\label{prop:mclean}
Infinitesimal special Lagrangian deformations of $L$ are in one to one
correspondence with cohomology classes in $H^1(L,\R)$. Moreover, the
deformations are unobstructed.
\end{proposition}

Specifically, a section of the normal bundle $v\in C^\infty(NL)$
determines a 1-form $\alpha=-\iota_v\omega\in\Omega^1(L,\R)$ and an
$(n-1)$-form $\beta=\iota_v\mathrm{Im}\,\Omega\in\Omega^{n-1}(L,\R)$.
These satisfy $\beta=\psi\,*_g\alpha$, where $\psi\in C^\infty(L,\R_+)$
is the ratio between the volume elements determined by $\Omega$ and $g$,
i.e.\ the norm of $\Omega$ with respect to the K\"ahler metric; moreover, the
deformation is
special Lagrangian if and only if $\alpha$ and $\beta$ are both closed.
Thus special Lagrangian deformations correspond to ``$\psi$-harmonic''
1-forms $$-\iota_v\omega\in \mathcal{H}^1_\psi(L)=\{\alpha\in
\Omega^1(L,\R)\,|\ 
d\alpha=0,\ d^*(\psi\alpha)=0\}.$$
In particular, special Lagrangian tori occur in $n$-dimensional families,
giving a local fibration structure provided that nontrivial $\psi$-harmonic
1-forms have no zeroes.

The base $B$ of a special Lagrangian torus fibration carries two natural
affine structures, which we call ``symplectic'' and ``complex''.
The first one, which encodes the symplectic geometry of $X$,
locally identifies $B$ with a domain
in $H^1(L,\R)$ ($L\approx T^n$).
At the level of tangent spaces, the cohomology class of 
$-\iota_v\omega$ provides an identification of $TB$ with
$H^1(L,\R)$; integrating, the local affine coordinates on $B$ are the
symplectic
areas swept by loops forming a basis of $H_1(L)$. The other affine
structure encodes the complex geometry of $X$, and locally identifies
$B$ with a domain in $H^{n-1}(L,\R)$. Namely, one uses the cohomology
class of $\iota_v\mathrm{Im}\,\Omega$ to identify
$TB$ with $H^{n-1}(L,\R)$, and the affine coordinates are obtained by
integrating $\mathrm{Im}\,\Omega$ over the $n$-chains swept by cycles
forming a basis of $H_{n-1}(L)$.

The dual special Lagrangian fibration can be constructed as a
moduli space $M$ of pairs $(L,\nabla)$, where $L\subset X$ is a special 
Lagrangian fiber and $\nabla$ is a rank 1 unitary local system over $L$.
The local geometry of $M$ is well-understood (cf.\ e.g.\ \cite{hitchin, leung,
gross}); in particular we have the following result (cf.\ e.g.\ \S 2 of \cite{Au}):

\begin{proposition}\label{prop:mirrgeom}
Let $M$ be the moduli space of pairs $(L,\nabla)$, where $L$ is a special
Lagrangian torus in $X$ and $\nabla$ is a flat $U(1)$ connection on
the trivial complex line bundle over $L$ up to gauge. Then
$M$ carries a natural integrable complex structure $J^\vee$ arising from
the identification $$T_{(L,\nabla)}M=\{(v,\alpha)\in C^\infty(NL)\oplus
\Omega^1(L,\R)\,|\,-\iota_v\omega+i\alpha\in \mathcal{H}^1_\psi(L)\otimes
\C\},$$ a holomorphic $n$-form
$$\Omega^\vee((v_1,\alpha_1),\dots,(v_n,\alpha_n))=\int_L
(-\iota_{v_1}\omega+i\alpha_1)\wedge\dots\wedge
(-\iota_{v_n}\omega+i\alpha_n),$$
and a compatible K\"ahler form
$$\omega^\vee((v_1,\alpha_1),(v_2,\alpha_2))=\int_L \alpha_2\wedge
\iota_{v_1}\Im\Omega-\alpha_1\wedge \iota_{v_2}\Im\Omega$$
(this formula for $\omega^\vee$ assumes that
$\int_L \Re\Omega$ has been suitably normalized).
\end{proposition}

In particular, $M$ can be viewed as a complexification of the moduli space
of special Lagrangian submanifolds; forgetting the connection gives a
projection map $f^\vee$ from $M$ to the real moduli space $B$. The fibers of
this projection are easily checked to be special Lagrangian tori in
the almost Calabi-Yau manifold $(M,J^\vee,\omega^\vee,\Omega^\vee)$. 

This special Lagrangian fibration on $M$ is fiberwise dual to the one
previously considered on $X$; they have the same base $B$, and passing 
from one fibration to the other simply amounts to exchanging the roles
of the two affine structures on $B$.
\medskip

In real life, unless we restrict ourselves to complex tori, we have to 
consider special Lagrangian torus fibrations with {\it singularities}. 
The base of the fibration is then a {\it singular} affine manifold, and 
the picture discussed above only holds away from the singularities. 
A natural idea would be to obtain the mirror by first constructing the
dual fibration away from the singularities, and then trying
to extend it over the singular locus. Unfortunately, this cannot be done
directly; instead we need to modify the complex geometry of $M$
by introducing {\it instanton corrections}.

To give some insight into the geometric meaning
of these corrections, consider the SYZ conjecture from the
perspective of homological mirror symmetry.
\medskip

Recall that Kontsevich's homological mirror symmetry conjecture
\cite{KoICM} predicts that the derived category of coherent sheaves
$D^b\mathrm{Coh}(X^\vee)$ of the mirror $X^\vee$ is equivalent to
the derived Fukaya category of $X$. For any point
$p\in X^\vee$, the skyscraper sheaf $\mathcal{O}_p$ is an object
of the derived category. Since $\mathrm{Ext}^*(\mathcal{O}_p,\mathcal{O}_p)
\simeq H^*(T^n;\C)$ (as a graded vector space), we expect that
$\mathcal{O}_p$ corresponds
to some object $\mathcal{L}_p$ of the derived Fukaya category of $X$
such that $\mathrm{End}(\mathcal{L}_p)\simeq H^*(T^n)$. It is
natural to conjecture that, generically, the object
$\mathcal{L}_p$ is a Lagrangian torus in $X$ with trivial
Maslov class, equipped with a rank 1 unitary local system,
and such that $HF^*(\mathcal{L}_p,\mathcal{L}_p)\simeq
H^*(T^n)$ (as a graded vector space). This suggests constructing
the mirror $X^\vee$ as a moduli space of such objects of the Fukaya 
category of $X$. (However it could still be the case that some points
of $X^\vee$ cannot be realized by honest Lagrangian tori in $X$.)

In the Calabi-Yau setting, it is expected that ``generically'' (i.e.,
subject to a certain stability condition) the Hamiltonian isotopy class
of the Lagrangian torus $\mathcal{L}_p$ should contain a
unique special Lagrangian representative \cite{thomas,thomas-yau}.
Hence it is natural to restrict one's attention to {\it special} 
Lagrangians, whose geometry is richer than that of Lagrangians: for instance, the moduli space 
considered in Proposition \ref{prop:mirrgeom} carries not only a complex
structure, but also a symplectic structure. However, if we only care about
the complex geometry of the mirror $X^\vee$ and not its symplectic geometry,
then it should not be necessary to consider special Lagrangians.

On the other hand, due to
wall-crossing phenomena, the ``convergent power series'' version of 
Lagrangian Floer homology which is directly relevant to the situation
here is not quite invariant under Hamiltonian isotopies (see e.g.\
\cite{cho-nonunitary}, and Section \ref{s:floer} below). Hence, we 
need to consider a {\it corrected}
equivalence relation on the moduli space of Lagrangian tori in $X$
equipped with unitary local systems. Loosely speaking, we'd like to say
that two Lagrangian tori (equipped unitary local systems) are equivalent
if they behave interchangeably with respect to convergent power series
Floer homology; however, giving a precise meaning to this statement is
rather tricky.

\subsection{Beyond the Calabi-Yau case: Landau-Ginzburg models}\label{ss:lg}
Assume now that $(X,J,\omega)$ is a K\"ahler manifold of complex dimension
$n$, and that $D\subset X$ is an effective divisor representing the
anticanonical class, with at most normal crossing singularities. Then
the inverse of the defining section of $D$ is a section of the canonical
line bundle $K_X$ over $X\setminus D$, i.e.\ a holomorphic volume form
$\Omega\in \Omega^{n,0}(X\setminus D)$ with simple poles along $D$.

We can try to construct a mirror to the almost Calabi-Yau manifold
$X\setminus D$ just as above, by considering a suitable moduli space 
of (special) Lagrangian tori in $X\setminus D$ equipped with unitary 
local systems.
The assumption on the behavior of $\Omega$ near $D$ is necessary for 
the existence of a special Lagrangian torus fibration with the
desired properties: for instance, 
a neighborhood of the origin in $\C$ equipped with
$\Omega=z^k\,dz$ does not contain any compact special Lagrangians
unless $k=-1$.

Compared to $X\setminus D$, the manifold $X$ contains essentially the
same Lagrangians. However, 
(special) Lagrangian tori in $X\setminus D$ typically
bound families of holomorphic discs in $X$, which causes their
Floer homology to be {\it obstructed}\/ in the
sense of Fukaya-Oh-Ohta-Ono \cite{FO3book}. Namely, Floer theory associates
to $\mathcal{L}=(L,\nabla)$ (where $L$ is a
Lagrangian torus in $X\setminus D$ and $\nabla$ is a flat $U(1)$ connection
on the trivial line bundle over $L$)
an element $\mathfrak{m}_0(\mathcal{L})\in CF^*(\mathcal{L},\mathcal{L})$, 
given by a weighted count of holomorphic discs in $(X,L)$.

More precisely,
recall that in Fukaya-Oh-Ohta-Ono's approach the Floer complex
$CF^*(\mathcal{L},\mathcal{L})$ is generated
by chains on $L$ (with suitable coefficients), and its element
$\mathfrak{m}_0(\mathcal{L})$ is defined as follows (see \cite{FO3book}
for details). 
Given a class $\beta\in \pi_2(X,L)$, the moduli space
$\mathcal{M}_k(L,\beta)$ of holomorphic discs in $(X,L)$ with $k$ boundary
marked points representing the class $\beta$ has expected dimension 
$n-3+k+\mu(\beta)$, where $\mu(\beta)$ is the Maslov index; when $L\subset
X\setminus D$ is special Lagrangian, $\mu(\beta)$ is simply twice the algebraic intersection number
$\beta\cdot [D]$ (see e.g.\ Lemma 3.1 of \cite{Au}). This moduli space can
be compactified by adding bubbled configurations. Assuming regularity, 
this yields a manifold with boundary, which carries a fundamental chain 
$[\overline{\mathcal{M}}_k(L,\beta)]$; otherwise, various
techniques can be used to define a virtual fundamental chain
$[\overline{\mathcal{M}}_k(L,\beta)]^{vir}$, usually dependent on auxiliary
perturbation data. The (virtual) fundamental chain of
$\overline{\mathcal{M}}_1(L,\beta)$ can be pushed forward by 
the evaluation map at the marked point, $ev:\overline{\mathcal{M}}_1(L,\beta)\to L$,
to obtain a chain in $L$: then one sets
\begin{equation}\label{eq:defm0}
\mathfrak{m}_0(\mathcal{L})=\sum_{\beta\in \pi_2(X,L)} z_\beta(\mathcal{L})\,
ev_*[\overline{\mathcal{M}}_1(L,\beta)]^{vir},\end{equation}
where the coefficient $z_\beta(\mathcal{L})$ reflects weighting by
symplectic area:
\begin{equation}\label{eq:zbeta}
z_\beta(\mathcal{L})=\exp(-\textstyle\int_\beta\omega)\,\mathrm{hol}_\nabla(\partial
\beta)\in\C^*,
\end{equation}
or $z_\beta(\mathcal{L})=T^{\int_\beta\omega}\,\mathrm{hol}_\nabla(\partial
\beta)\in\Lambda_0$
if using Novikov coefficients to avoid convergence issues 
(see below).

Note that $z_\beta$ as defined by (\ref{eq:zbeta})
is locally a  holomorphic function with respect to the complex structure
$J^\vee$ introduced in Proposition \ref{prop:mirrgeom}. Indeed, recall
that the tangent space to the moduli space $M$ is identified with the space
of complex-valued $\psi$-harmonic 1-forms on $L$; the differential of $\log
z_\beta$ is just the linear form on $\mathcal{H}^1_\psi(L)\otimes\C$
given by integration on
the homology class $\partial\beta\in H_1(L)$.

In this paper we will mostly consider {\it weakly unobstructed}\/ Lagrangians,
i.e.\ those for which $\mathfrak{m}_0(\mathcal{L})$ is a multiple of the unit 
(the fundamental cycle of $L$). In that case, the Floer differential on
$CF^*(\mathcal{L},\mathcal{L})$ does square to zero, but given two
Lagrangians $\mathcal{L},\mathcal{L}'$ we find that $CF^*(\L,\L')$ may
not be well-defined as a chain complex. To understand the obstruction,
recall that the count of holomorphic triangles equips $CF^*(\L,\L')$ with
the structure of a left module over $CF^*(\L,\L)$ and a right module over
$CF^*(\L',\L')$. Writing $\m_2$ for both module maps, an analysis of the
boundary of 1-dimensional moduli spaces shows that the differential
on $CF^*(\mathcal{L},\mathcal{L}')$ squares to 
$$\m_2(\m_0(\L'),\cdot)-\m_2(\cdot,\m_0(\L)).$$ The assumption that
$\m_0$ is a multiple of the identity implies that
Floer homology is only defined for pairs of Lagrangians
which have the same $\mathfrak{m}_0$. Moreover, even though
the Floer homology group $HF^*(\mathcal{L},\mathcal{L})$ can still be
defined, it is generically zero due to contributions of holomorphic 
discs in $(X,L)$ to the Floer differential; in that case $\mathcal{L}$ 
is a trivial object of the Fukaya category. 

On the mirror side, these features of the theory
can be replicated by the introduction of a {\it superpotential}, i.e.\ a
holomorphic function $W:X^\vee\to\C$ on the mirror of $X\setminus D$. 
$W$ can be thought of as an obstruction term for the B-model on $X^\vee$,
playing the same role as $\mathfrak{m}_0$ for the A-model on $X$. More
precisely, homological mirror symmetry predicts that the derived Fukaya
category of $X$ is equivalent to the {\it derived category of singularities}
of the mirror Landau-Ginzburg model $(X^\vee,W)$ \cite{KL,orlov}. 
This category is actually
a {\it collection} of categories indexed by complex numbers, just as
the derived Fukaya category of $X$ is a collection of categories indexed
by the values of $\mathfrak{m}_0$. 

Given $\lambda\in\C$, one defines
$D^b_{sing}(W,\lambda)=D^b\mathrm{Coh}(W^{-1}(\lambda))/
\mathrm{Perf}(W^{-1}(\lambda))$, the quotient of the derived category of
coherent sheaves on the fiber $W^{-1}(\lambda)$ by the subcategory of
perfect complexes. Since for smooth fibers the derived category
of coherent sheaves is generated by vector bundles, this quotient is trivial
unless $\lambda$ is a critical value of $W$;
in particular, a point of $X^\vee$ defines a nontrivial object of the
derived category of singularities only if it is a critical point of $W$.
Alternatively, this category can also be defined in terms of {\it matrix
factorizations}. Assuming $X^\vee$ to be affine for simplicity, a matrix
factorization is a $\Z/2$-graded projective $\C[X^\vee]$-module together
with an odd endomorphism $\delta$ such that $\delta^2=(W-\lambda)\,\mathrm{id}$.
For a fixed value of $\lambda$, matrix factorizations yield a $\Z/2$-graded 
dg-category, whose cohomological category is equivalent to
$D^b_{sing}(W,\lambda)$ by a result of Orlov \cite{orlov}. However,
if we consider two matrix factorizations $(P_1,\delta_1)$ and
$(P_2,\delta_2)$ associated to two values $\lambda_1,\lambda_2\in\C$, 
then the differential on $\mathrm{hom}((P_1,\delta_1),(P_2,\delta_2))$
squares to $(\lambda_1-\lambda_2)\,\mathrm{id}$, similarly to the Floer differential
on the Floer complex of two Lagrangians with different values of
$\mathfrak{m}_0$. 

This motivates the following conjecture:

\begin{conj}\label{conj:mirror}
The mirror of $X$ is the {\it Landau-Ginzburg model} $(X^\vee,W)$,
where
\begin{enumerate}
\item $X^\vee$ is a mirror of the almost Calabi-Yau manifold 
$X\setminus D$, i.e.\ a (corrected and completed) moduli space of special Lagrangian
tori in $X\setminus D$ equipped with rank 1 unitary local systems;
\medskip
\item $W:X^\vee\to\C$ is a holomorphic function defined as follows:
if $p\in X^\vee$ corresponds to a special Lagrangian $\mathcal{L}_p=
(L,\nabla)$, then
\begin{equation}\label{eq:W}
W(p)=\sum_{\beta\in\pi_2(X,L),\ \mu(\beta)=2}\!\!\!
n_\beta(\mathcal{L}_p)\,z_\beta(\mathcal{L}_p),\end{equation}
where $n_\beta(\mathcal{L}_p)$ is the degree of the evaluation chain
$ev_*[\overline{\mathcal{M}}_1(L,\beta)]^{vir}$, i.e., the (virtual) number of
holomorphic discs in the class $\beta$ passing through a generic point
of $L$, and the weight $z_\beta(\mathcal{L}_p)$ is given by (\ref{eq:zbeta}).
\end{enumerate}
\end{conj}

There are several issues with the formula (\ref{eq:W}). To start with,
except in specific cases (e.g.\ Fano toric varieties), there is no
guarantee that the sum in (\ref{eq:W}) converges. The rigorous way to
deal with this issue is to work over the Novikov ring 
\begin{equation}\label{eq:novikov}
\Lambda_0=\left\{\textstyle\sum\limits_i\, a_i\, T^{\lambda_i}\,\Big|\,
a_i\in\C,\ \lambda_i\in\R_{\ge 0},\ \lambda_i\to +\infty\right\}
\end{equation}
rather than over complex numbers.
Holomorphic discs in a class $\beta$ are then counted 
with weight $T^{\int_\beta\omega}\,\mathrm{hol}_\nabla(\partial\beta)$
instead of $\exp(-\int_\beta\omega)\,\mathrm{hol}_\nabla(\partial\beta)$.
Assuming convergence, setting $T=e^{-1}$ recovers the complex coefficient
version.

Morally, working over Novikov coefficients simultaneously encodes
the family of mirrors for $X$ equipped with the family of K\"ahler forms
$\kappa\omega,\ \kappa\in\R_+$. Namely, the mirror manifold should be
constructed as a variety defined over the Novikov field $\Lambda$ (the
field of fractions of $\Lambda_0$), and the 
superpotential as a regular function with values in $\Lambda$.
If convergence holds, then setting $T=\exp(-\kappa)$ recovers the
complex mirror to $(X,\kappa\omega)$; if convergence fails for all
values of $T$, the mirror might actually exist only in a formal sense
near the {\it large volume limit\/} $\kappa\to\infty$.

Another issue with Conjecture \ref{conj:mirror} is the definition of
the numbers $n_\beta(\mathcal{L}_p)$. Roughly speaking, the value of
the superpotential is meant to be ``the coefficient of the fundamental
chain $[L]$ in $\mathfrak{m}_0$''. However, in real life not all Lagrangians
are weakly unobstructed: due to bubbling of Maslov
index 0 discs, for a given class $\beta$ with $\mu(\beta)=2$ the
chain $ev_*[\overline{\mathcal{M}}_1(L,\beta)]^{vir}$ is 
in general not a cycle. Thus we can still define $n_\beta$ to be its ``degree'' (or
multiplicity) at some point $q\in L$, but the answer depends on the choice 
of $q$. Alternatively, we can complete the chain to a cycle, e.g.\
by choosing a ``weak bounding cochain'' in the sense of
Fukaya-Oh-Ohta-Ono \cite{FO3book}, or more geometrically, by 
considering not only holomorphic discs but also holomorphic ``clusters''
in the sense of Cornea-Lalonde \cite{CL}; however, $n_\beta$ will then
depend on some auxiliary data (in the cluster approach, a Morse function 
on $L$). 

Even if we equip each $L$ with the appropriate auxiliary data (e.g.\ a
base point or a Morse function),
the numbers $n_\beta$ will typically vary in a discontinuous manner due to
wall-crossing phenomena. However, recall that $X^\vee$ differs from the
naive moduli space of Proposition \ref{prop:mirrgeom} 
by instanton corrections. Namely, $X^\vee$ is more 
accurately described as a (completed) moduli space of Lagrangian tori 
$L\subset X\setminus D$ equipped
with not only a $U(1)$ local system but also the auxiliary data needed to
make sense of the Floer theory of $L$ in general and of the numbers
$n_\beta$ in particular. The equivalence relation on this set of Lagrangians
equipped with extra data is Floer-theoretic in nature. General considerations about wall-crossing
and continuation maps in Floer theory imply that, even though
the individual numbers $n_\beta$ depend on the choice of
a representative in the equivalence class, by construction the 
superpotential $W$ given
by (\ref{eq:W}) is a single-valued smooth function on the corrected moduli
space. The reader is referred to \S 19.1 in \cite{FO3book} and \S 3 
in \cite{Au} for details.

In this paper we will assume that things don't go {\it completely}\ wrong,
namely that our Lagrangians are weakly unobstructed except when they
lie near a certain collection of walls in the moduli space. In this case,
the process which yields the corrected moduli space from
the naive one can be thought of decomposing $M$ into {\it chambers} over
which the $n_\beta$ are locally constant, and gluing these chambers by
{\it analytic} changes of coordinates dictated by the enumerative geometry
of Maslov index 0 discs on the wall. Thus, the analyticity of $W$ on the
corrected mirror follows from that of $z_\beta$ on the uncorrected moduli
space.

One last thing to mention is that the incompleteness of the K\"ahler
metric on $X\setminus D$ causes the moduli space of Lagrangians to be
similarly incomplete. This is readily apparent if we observe that,
since $|z_\beta|=\exp(-\int_\beta\omega)$, each variable $z_\beta$
appearing in the sum (\ref{eq:W}) takes values in the unit disc.
We will want to define the mirror of $X$ to be a
{\it larger} space, obtained by analytic continuation of the
instanton-corrected moduli space of Lagrangian tori (i.e., roughly
speaking, allowing $|z_\beta|$ to be arbitrarily large).
One way to think of the points of $X^\vee$ added in the completion
process is as Lagrangian tori in $X\setminus D$ equipped with
{\it non-unitary} local systems; however this can lead to serious 
convergence issues, even when working over the Novikov ring.

Another way to think about the completion process, under the
assumption that $D$ is nef, is in terms of {\it inflating} $X$ along
$D$, i.e\ replacing the K\"ahler form $\omega$ by $\omega_t=\omega+t\eta$
where the $(1,1)$-form $\eta$ is Poincar\'e dual to $D$ and supported in
a neighborhood of $D$; this ``enlarges'' the moduli space of Lagrangians
near $D$, and simultaneously increases the area of all Maslov index 2 
discs by $t$, i.e.\ rescales the superpotential by a factor of $e^{-t}$.
Taking the limit as $t\to\infty$ (and rescaling the superpotential
appropriately) yields the completed mirror.

\subsection{Example: Fano toric varieties} \label{ss:toric}

Let $(X,\omega,J)$ be a smooth toric variety of complex dimension $n$.
In this section we additionally assume that $X$ is Fano, i.e.\ its
anticanonical divisor is ample.
As a K\"ahler manifold, $X$ is determined by its moment
polytope $\Delta\subset\R^n$, a convex polytope in which every
facet admits an integer normal vector, $n$ facets meet at
every vertex, and their primitive integer normal vectors form a basis of
$\Z^n$. The moment map $\phi:X\to\R^n$ identifies the orbit space of the
$T^n$-action on $X$ with $\Delta$. From the point of view of complex
geometry, the preimage of the interior of $\Delta$ is an open dense
subset $U$ of $X$, biholomorphic to $(\C^*)^n$, on which $T^n=(S^1)^n$ acts
in the standard manner. Moreover $X$ admits an open cover by affine
subsets biholomorphic to $\C^n$, which are the
preimages of the open stars of the vertices of $\Delta$ (i.e., the union
of all the strata whose closure contains the given vertex).

For each facet $F$ of $\Delta$, the preimage $\phi^{-1}(F)=D_F$ is a
hypersurface in $X$; the union of these hypersurfaces defines the toric
anticanonical divisor $D=\sum_F D_F$. The standard holomorphic volume 
form on $(\C^*)^n\simeq U=X\setminus D$, defined in coordinates by
$\Omega=d\log x_1\wedge\dots\wedge d\log x_n$, determines a section of
$K_X$ with poles along $D$.

It is straightforward to check that the orbits of the 
$T^n$-action are special Lagrangian with respect to $\omega$ and $\Omega$. Thus the moment map determines a special
Lagrangian fibration on $X\setminus D$, with base $B=\mathrm{int}(\Delta)$; by
definition, the symplectic affine structure induced 
on $B$ by the identification $TB\approx H^1(L,\R)$ is precisely the
standard one coming from the inclusion of $B$ in $\R^n$ (up to
a scaling factor of $2\pi$).

Consider a $T^n$-orbit $L$ in the open stratum $X\setminus D\approx
(\C^*)^n$, and a flat $U(1)$-connection $\nabla$ on the trivial bundle
over $L$. Let
$$z_j(L,\nabla)=\exp(-2\pi\phi_j(L))\,\mathrm{hol}_\nabla(\gamma_j),$$
where $\phi_j$ is the $j$-th component of the moment map, i.e.\ the
Hamiltonian for the action of the $j$-th factor of $T^n$, and
$\gamma_j=[S^1(r_j)]\in H_1(L)$ is the homology class corresponding
to the $j$-th factor in $L=S^1(r_1)\times\dots\times S^1(r_n)\subset
(\C^*)^n$. Then $z_1,\dots,z_n$ are holomorphic coordinates on the moduli
space $M$ of pairs $(L,\nabla)$ equipped with the complex structure
$J^\vee$ of Proposition \ref{prop:mirrgeom}.

For each facet $F$ of $\Delta$, denote by $\nu(F)\in\Z^n$ the primitive
integer normal vector to $F$ pointing into $\Delta$, and let
$\alpha(F)\in\R$
be the constant such that the equation of $F$ is $\langle
\nu(F),\phi\rangle+
\alpha(F)=0$. Moreover, given $a=(a_1,\dots,a_n)\in \Z^n$ we denote by
$z^a$ the Laurent monomial $z_1^{a_1}\dots z_n^{a_n}$.

\begin{proposition}\label{prop:toric}
The SYZ mirror to the smooth Fano toric variety $X$ is\/ $(\C^*)^n$ equipped
with a superpotential given by the Laurent polynomial
\begin{equation}\label{eq:toricW}
W=\sum_{F\ \mathrm{facet}} e^{-2\pi\alpha(F)}\,z^{\nu(F)}.
\end{equation}
\end{proposition}

\noindent
More precisely, the moduli space $M$ of pairs $(L,\nabla)$ is biholomorphic
to the bounded open subset of $(\C^*)^n$ consisting of all points 
$(z_1,\dots,z_n)$ such that each term in the sum $(\ref{eq:toricW})$ has 
norm less than 1; however, the completed mirror is all of $(\C^*)^n$.

Proposition \ref{prop:toric} is a well-known result, which appears in
many places; for completeness we give a very brief sketch of a geometric
proof (see also \cite{hori,cho-oh,Au,FO3toric} for more details).

\proof[Sketch of proof]
Consider a pair $(L,\nabla)$ as above, and recall that $L$ can be identified
with a product torus $S^1(r_1)\times \dots\times S^1(r_n)$ inside $(\C^*)^n$.
It follows from the maximum principle that $L$ does not bound any
non-constant holomorphic discs in $(\C^*)^n$; since the Maslov index of 
a disc in $(X,L)$ is twice its intersection number with the toric divisor
$D$, this eliminates the possibility of Maslov index 0 discs.
Moreover, since $X$ is Fano, all holomorphic spheres
in $X$ have positive Chern number. It follows that the moduli spaces of
Maslov index 2 holomorphic discs in $(X,L)$ are all compact,
and that we do not have to worry about possible contributions from 
bubble trees of total Maslov index 2; this is in sharp contrast with 
the non-Fano case, see \S \ref{ss:hirzebruch}.

A holomorphic disc of Maslov index 2 in $(X,L)$ intersects $D$ at a single
point, and in
particular it intersects only one of the components, say $D_F$ for some
facet $F$ of $\Delta$. 
Cho and Oh \cite{cho-oh} observed that for each facet $F$ there is a unique such disc
whose boundary passes through a given point $\mathbf{x}^0=
(x^0_1,\dots,x^0_n)\in L$; in terms of the components
$(\nu_1,\dots,\nu_n)$ of the normal vector $\nu(F)$, this disc can be
parametrized by the map 
\begin{equation}\label{eq:holdisc}
w\mapsto (w^{\nu_1} x^0_1,\dots,w^{\nu_n} x^0_n)\end{equation}
(for $w\in D^2\setminus\{0\}$; the point $w=0$ corresponds to the
intersection with $D_F$).

This is easiest to check in the model case where $X=\C^n$, the
moment polytope is the positive octant $\R_{\ge 0}^n$, and the normal
vectors to the facets form the standard basis of $\Z^n$. Namely, the maximum
principle implies that holomorphic discs of Maslov index 2 with boundary 
in a product torus in $\C^n$ are given by
maps with only one non-constant component, and up to reparametrization
that non-constant component can be assumed to be linear.
The general case is proved by working in an affine chart centered
at a vertex of $\Delta$ adjacent to the considered facet $F$, and using
a suitable change of coordinates to reduce to the previous case.

A careful calculation shows that the map (\ref{eq:holdisc}) is regular,
and that its contribution to the signed count of holomorphic discs is $+1$.
Moreover, it follows from the definition of the moment map that the 
symplectic area of this disc is $2\pi (\langle \nu(F),\phi(L)\rangle+\alpha(F))$.
(This is again easiest to check in the model case of $\C^n$; the general
case follows by performing a suitable change of coodinates). Exponentiating
and multiplying by the appropriate holonomy factor, one arrives at 
(\ref{eq:toricW}).

Finally, recall that the interior of $\Delta$ is defined by the inequalities
$\langle \nu(F),\phi(L)\rangle+\alpha(F)>0$ for every facet $F$;
exponentiating, this corresponds exactly to the constraint that
$|e^{-2\pi\alpha(F)}z^{\nu(F)}|<1$ for every facet $F$. However,
adding the Poincar\'e dual of $tD$ to $\omega$ enlarges the moment
polytope by $t/2\pi$ in every direction, i.e.\ it increases $\alpha(F)$ by 
$t/2\pi$ for all facets. This makes $M$ a larger subset of $(\C^*)^n$;
rescaling $W$ by $e^t$ and taking the limit as $t\to +\infty$, we obtain
all of $(\C^*)^n$ as claimed.
\endproof

\section{Examples of wall-crossing and instanton corrections} \label{s:examples}

\subsection{First examples}\label{ss:examples}
In this section we give two simple examples illustrating the construction
of the mirror and the process of instanton corrections. The first example
is explained in detail in \S 5 of \cite{Au}, while the second example is the
starting point of~\cite{AAK}; the two examples are in fact very similar.

\subsubsection{}\label{ex:conic}

Consider $X=\C^2$, equipped with a toric K\"ahler form $\omega$ and
the holomorphic volume form $\Omega=dx\wedge dy/(xy-\epsilon)$, which
has
poles along the conic $D=\{xy=\epsilon\}$. Then $X\setminus D$ carries
a fibration by special Lagrangian tori $$T_{r,\lambda}=\{(x,y)\in\C^2,\ 
|xy-\epsilon|=r,\ \mu_{S^1}(x,y)=\lambda\},$$
where $\mu_{S^1}$ is the moment map for the $S^1$-action $e^{i\theta}\cdot
(x,y)=(e^{i\theta}x,e^{-i\theta}y)$, for instance
$\mu_{S^1}(x,y)=\frac12(|x|^2-|y|^2)$ for $\omega=\frac{i}{2}(dx\wedge
d\bar{x}+dy\wedge d\bar{y})$.
These tori are most easily visualized in terms of the
projection $f:(x,y)\mapsto xy$, whose fibers are affine conics,
each of which carries a $S^1$-action.
The torus $T_{r,\lambda}$ lies in the preimage by $f$ of a circle of radius $r$
centered at $\epsilon$, and consists of a single $S^1$-orbit inside each
fiber. In particular, $T_{r,\lambda}$ is smooth unless
$(r,\lambda)=(|\epsilon|,0)$, where we have a nodal singularity at the
origin.

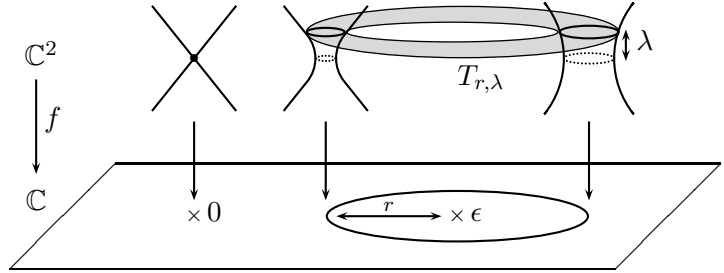
\begin{figure}[h]
\setlength{\unitlength}{7mm}
\begin{picture}(13.5,5)(0,0.2)
\psset{unit=\unitlength}
\newgray{ltgray}{0.85}
\psellipse[linewidth=0.5pt,fillstyle=solid,fillcolor=ltgray](8.6,4.5)(2.98,0.5)
\psellipse[linewidth=0.5pt,fillstyle=solid,fillcolor=white](8.4,4.5)(2.02,0.2)
\put(0,0){\line(1,0){11.5}}
\put(0,0){\line(1,1){2}}
\put(2,2){\line(1,0){11.5}}
\put(11.5,0){\line(1,1){2}}
\put(3.5,1){\makebox(0,0)[cc]{\tiny$\times$}}
\put(3.75,0.9){\small $0$}
\put(8.5,1){\makebox(0,0)[cc]{\tiny$\times$}}
\put(8.75,0.9){\small $\epsilon$}
\psline{<->}(6.2,1)(8.2,1)
\put(7.1,1.1){\tiny $r$}
\psellipse(8.5,1)(2.5,0.5)
\psline[linearc=0](2.7,3)(3.5,4)(2.7,5)
\psline[linearc=0](4.3,3)(3.5,4)(4.3,5)
\psline[linearc=0.7](5.2,3)(6,4)(5.2,5)
\psline[linearc=0.7](6.8,3)(6,4)(6.8,5)
\psline[linearc=1.7](10.2,3)(11,4)(10.2,5)
\psline[linearc=1.7](11.8,3)(11,4)(11.8,5)
\put(3.5,4){\circle*{0.15}}
\psellipse[linestyle=dotted,dotsep=0.6pt](6,4)(0.2,0.07)
\psellipse[linestyle=dotted,dotsep=0.6pt](11,4)(0.5,0.12)
\psellipse(11,4.5)(0.58,0.14)
\psellipse(6,4.5)(0.38,0.1)
\psline{<->}(11.7,3.95)(11.7,4.55)
\put(11.9,4.15){\small $\lambda$}
\psline{->}(3.5,2.8)(3.5,1.3)
\psline{->}(6,2.8)(6,1.3)
\psline{->}(11,2.8)(11,1.3)
\put(0.5,1.3){\makebox(0,0)[cc]{\small $\C$}}
\put(0.5,4.1){\makebox(0,0)[cc]{\small $\C^2$\!\!}}
\put(0.65,2.7){\small $f$}
\psline{->}(0.5,3.6)(0.5,1.8)
\put(8.5,3.5){\small $T_{r,\lambda}$}
\end{picture}
\caption{The special Lagrangian torus $T_{r,\lambda}$ in $\C^2\setminus D$}
\label{fig:conic}
\end{figure}

\noindent
One can check that $T_{r,\lambda}$ is special Lagrangian either by direct
calculation, or by observing that $T_{r,\lambda}$
is the lift of a special Lagrangian circle in the reduced space 
$X_{red,\lambda}=\mu^{-1}_{S^1}(\lambda)/S^1$
equipped with the reduced K\"ahler form $\omega_{red,\lambda}$ and the
reduced holomorphic volume form $\Omega_{red,\lambda}=
\iota_{(\partial/\partial\theta)^\#}\Omega=i\,d\log(xy-\epsilon)$;
see \S 5 of \cite{Au}.

As seen in \S \ref{s:background}, away from $(r,\lambda)=(|\epsilon|,0)$
the
moduli space $M$ of pairs consisting of a torus
$L=T_{r,\lambda}$ and a $U(1)$ local system $\nabla$ carries a natural
complex structure, for which the functions
$z_\beta=\exp(-\int_\beta\omega)\,\hol_\nabla(\partial\beta)$,
$\beta\in\pi_2(X,L)$ are holomorphic.

Wall-crossing occurs at $r=|\epsilon|$, namely the tori
$T_{|\epsilon|,\lambda}$ for $\lambda>0$ intersect the $x$-axis in
a circle, and thus bound a holomorphic disc contained in the fiber
$f^{-1}(0)$, which has Maslov index 0. Denote by $\alpha$ the relative
homotopy class of this disc, and by $w=z_\alpha$ the corresponding
holomorphic weight, which satisfies $|w|=e^{-\lambda}$. 
Similarly the tori $T_{|\epsilon|,\lambda}$ for $\lambda<0$ bound
a holomorphic disc contained in the $y$-axis, representing the
relative class $-\alpha$ and with associated weight $z_{-\alpha}=w^{-1}$.

Since the projection $f$ is holomorphic, holomorphic discs of Maslov 
index 2 in $(\C^2,T_{r,\lambda})$ are sections of $f$ over the disc 
of radius $r$ centered at $\epsilon$. When $r>|\epsilon|$ there are 
two families of such discs; these can be found either by explicit
calculation, or by deforming $T_{r,\lambda}$ to a product
torus $S^1(r_1)\times S^1(r_2)$ (by deforming the circle centered at 
$\epsilon$ to a circle centered at the origin, without crossing $\epsilon$),
for which the discs are simply $D^2(r_1)\times \{y\}$ and $\{x\}\times
D^2(r_2)$. Denote by $\beta_1$ and $\beta_2$ respectively the classes
of these discs, and by $z_1$ and $z_2$ the corresponding weights,
which satisfy $z_1/z_2=w$. In terms of these coordinates on $M$ the
superpotential is then given by $W=z_1+z_2$.

On the other hand, when $r<|\epsilon|$ there is only one family of
Maslov index 2 discs in $(\C^2,T_{r,\lambda})$. This is easiest to see
by deforming $T_{r,\lambda}$ to the {\it Chekanov torus} $|xy-\epsilon|=r$,
$|x|=|y|$ (if $\omega$ is invariant under $x\leftrightarrow y$ this is 
simply $T_{r,0}$); then the maximum principle applied to $y/x$ implies that 
Maslov index 2 discs are portions of lines $y=ax$, $|a|=1$.
Denoting by $\beta_0$ the class of this disc, and by $u$ the corresponding
weight, in the region $r<|\epsilon|$ the superpotential is $W=u$.

When we increase the value of $r$ past $r=|\epsilon|$, for $\lambda>0$,
the family of holomorphic discs in the class $\beta_0$ deforms naturally 
into the family of discs in the class $\beta_2$ mentioned above; the
coordinates on $M$ naturally glue according to $u=z_2$, $w=z_1/z_2$.
On the other hand, for $\lambda<0$ the class $\beta_0$ deforms naturally
into the class $\beta_1$, so that the coordinates glue according to
$u=z_1$, $w=z_1/z_2$. The discrepancy in these gluings is due to
the monodromy of our special
Lagrangian fibration around the singular fiber $T_{|\epsilon|,0}$, which
acts nontrivially on $\pi_2(\C^2,T_{r,\lambda})$:
while the coordinate $w$ is defined globally on $M$, 
$z_1$ and $z_2$ do not extend to global coordinates.

There are now two issues: the complex manifold $M$ does not extend 
across the singularity at $(r,\lambda)=(|\epsilon|,0)$, and the
superpotential $W$ is discontinuous across the walls. Both issues are
fixed simultaneously by instanton corrections. Namely, we correct
the coordinate change across the wall $r=|\epsilon|$, $\lambda>0$ 
to $u=z_2(1+w)$. The correction factor $1+w$
indicates that, upon deforming $T_{r,\lambda}$ by increasing the value 
of $r$ past $|\epsilon|$, Maslov index 2 discs in the class $\beta_0$ 
give rise not only to discs in the class $\beta_2$ (by a straightforward
deformation), but also to new discs in the class $\beta_1=\beta_2+\alpha$ 
formed by attaching the exceptional disc bounded by
$T_{|\epsilon|,\lambda}$. Similarly, across $r=|\epsilon|$, $\lambda<0$,
we correct the gluing to $u=z_1(1+w^{-1})$, to take into account the
exceptional disc in the class $-\alpha$ bounded by $T_{|\epsilon|,\lambda}$.

The corrected gluings both come out to be $u=z_1+z_2$, which means that
we now have a well-defined mirror $X^\vee$, carrying a well-defined 
superpotential $W=u=z_1+z_2$. More precisely, using the coordinates $(u,w)$
on the chamber $r<|\epsilon|$, and the coordinates $(v,w)$ with $v=z_2^{-1}$
and $w=z_1/z_2$ on the chamber $r>|\epsilon|$, we claim that the corrected
and completed mirror is
$$X^\vee=\{(u,v,w)\in \C^2\times \C^*,\ uv=1+w\},\qquad W=u.$$
More precisely, the region $r>|\epsilon|$ of our special Lagrangian
fibration corresponds to $|z_1|$ and $|z_2|$ small, i.e.\ $|v|$ large; whereas
the region $r<|\epsilon|$ corresponds to $|u|$ large compared to
$e^{-|\epsilon|}$. When considering $M$ we also have $|u|<1$, as $|u|\to 1$
corresponds to $r\to 0$, but this constraint is
removed by the completion process, which enlarges $X$ along the conic
$xy=\epsilon$ by symplectic inflation. It turns out that we also have
to complete $X^\vee$ in the ``intermediate'' region where $u$ and $v$
are both small, in particular allowing these variables to vanish; for
otherwise, the corrected mirror would have ``gaps'' in the heavily corrected
region near $(r,\lambda)=(|\epsilon|,0)$. Let us also point out that $X^\vee$
is again the complement of a conic in $\C^2$.

General features of wall-crossing in Floer theory ensure that, when crossing
a wall, holomorphic disc counts (and hence the superpotential) can be made 
to match by introducing a suitable analytic change of coordinates,
consistently for all homotopy classes (see \S 19.1 of \cite{FO3book} and
\S 3 of \cite{Au}). For instance, if we compactified
$\C^2$ to $\CP^2$ or $\CP^1\times\CP^1$, then the tori $T_{r,\lambda}$
would bound additional families of Maslov index 2 holomorphic discs (passing
through the divisors at infinity), leading to additional terms in the
superpotential; however, these terms also match under the corrected gluing 
$u=z_1+z_2$ (see \S 5 of \cite{Au}).

\subsubsection{} \label{ex:c2blowup} Consider $\C^2$ equipped with 
the standard holomorphic volume form $d\log x\wedge d\log y$ (with poles along
the coordinate axes), and blow up the point $(1,0)$. This yields
a complex manifold $X$ equipped with the holomorphic volume form
$\Omega=\pi^*(d\log x\wedge d\log y)$, with poles along the proper
transform $D$ of the coordinate axes. Observe that the
$S^1$-action $e^{i\theta}\cdot(x,y)=
(x,e^{i\theta}y)$ lifts to $X$, and consider an $S^1$-invariant
K\"ahler form $\omega$ for which the area of the exceptional divisor is
$\epsilon$. Denote by
$\mu_{S^1}:X\to\R$ the moment map for the $S^1$-action, normalized to equal
0 on the proper transform of the $x$-axis and $\epsilon$ at the isolated fixed point.
Then the $S^1$-invariant tori $$L_{r,\lambda}=\{|\pi^*x|=r,\
\mu_{S^1}=\lambda\}$$
define a special Lagrangian fibration on $X\setminus D$,
with a nodal singularity at the isolated fixed point (for
$(r,\lambda)=(1,\epsilon)$) \cite{AAK}.

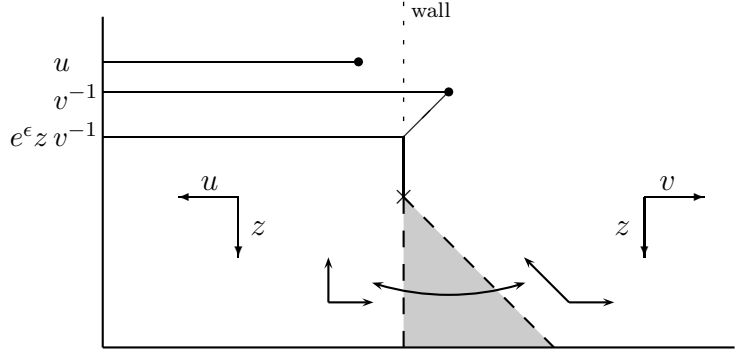
\begin{figure}[t]
\setlength{\unitlength}{2cm}
\begin{picture}(4.5,2.4)(-2,-2)
\psset{unit=\unitlength}
\newgray{ltgray}{0.8}
\pspolygon[fillstyle=solid,fillcolor=ltgray,linestyle=none](0,-1)(0,-2)(1,-2)
\psline(-2,-2)(2.2,-2)
\psline(-2,-2)(-2,0.2)
\put(0,-1){\makebox(0,0)[cc]{\small $\times$}}
\psline[linestyle=dashed,dash=0.1 0.1](0,-1)(0,-2)
\psline[linestyle=dashed,dash=0.1 0.1](0,-1)(1,-2)
\psline[linestyle=dashed,dash=0.03 0.1,linewidth=0.5pt](0,0.3)(0,-1)
\put(0.05,0.2){\tiny wall}
\psarc{<->}(0.3,0){1.65}{-108}{-72}
\psline{->}(-0.5,-1.7)(-0.5,-1.4)
\psline{->}(-0.5,-1.7)(-0.2,-1.7)
\psline{->}(1.1,-1.7)(1.4,-1.7)
\psline{->}(1.1,-1.7)(0.8,-1.4)
\put(-1.1,-1){\vector(-1,0){0.4}}
\put(-1.1,-1){\vector(0,-1){0.4}}
\put(-1.02,-1.25){$z$}
\put(-1.35,-0.95){$u$}
\put(1.6,-1){\vector(1,0){0.4}}
\put(1.6,-1){\vector(0,-1){0.4}}
\put(1.4,-1.25){$z$}
\put(1.7,-0.95){$v$}
\pscircle*(-0.3,-0.1){0.03}
\put(-0.3,-0.1){\line(-1,0){1.7}}
\multiput(-1.55,-0.1)(-0.1,0){4}{\line(-1,0){0.05}}
\put(-2.33,-0.16){\small $u$}
\pscircle*(0.3,-0.3){0.03}
\put(0.3,-0.3){\line(-1,-1){0.3}}
\put(0,-0.6){\line(0,-1){0.4}}
\put(0,-0.6){\line(-1,0){2}}
\put(0.3,-0.3){\line(-1,0){2.3}}
\put(-2.33,-0.41){\small $v^{-1}$}
\put(-2.62,-0.66){\small $e^{\epsilon}z\,v^{-1}$}
\end{picture}
\caption{A special Lagrangian fibration on the blowup of $\C^2$}
\label{fig:blowup}
\end{figure}

The base of this special Lagrangian fibration is pictured on
Figure \ref{fig:blowup}, where the vertical axis corresponds to the moment
map, and a cut has been made below the singular point to depict the 
monodromy of the symplectic affine structure.

For $r=1$ the Lagrangian tori $L_{r,\lambda}$ 
bound exceptional holomorphic discs, which causes wall-crossing:
for $\lambda>\epsilon$, $L_{1,\lambda}$ bounds a Maslov index 0 disc in
the proper transform of the line $x=1$, whereas for $\lambda<\epsilon$,
$L_{1,\lambda}$ splits the exceptional divisor of the blowup into two
discs, one of which has Maslov index 0. Thus,
we have to consider the chambers $r>1$ and $r<1$ separately.

When $r<1$, the Lagrangian torus $L_{r,\lambda}$ bounds two families
of Maslov index 2 discs. One family consists of the portions where
$\mu_{S^1}<\lambda$ of the lines $x=\mathrm{constant}$; we denote by
$\delta$ the homotopy class of these discs, and by $z$ ($=z_\delta$) the
corresponding holomorphic coordinate on $M$, which satisfies
$|z|=e^{-\lambda}$. The other family
consists of discs intersecting the $y$-axis, and is easiest to see by
deforming $L_{r,\lambda}$ to a product torus, upon which it becomes the
family of discs of radius $r$ in the lines $y=\mathrm{constant}$. 
(In fact, $L_{r,\lambda}$ is typically already a product torus for
$r$ sufficiently different from 1, when it lies in the region where
the blow-up operation does not affect the K\"ahler form.) We denote
by $\beta$ the class of these discs, and by $u$ the corresponding
holomorphic coordinate on $M$. The coordinates $u$ and $z$ on $M$ can
be thought of as (exponentiated) complexifications of the
affine coordinates on the base pictured on Figure \ref{fig:blowup}.

On the other hand, when $r>1$ the torus $L_{r,\lambda}$ bounds three
families of Maslov index 2 discs. As before, one of these families consists
of the portions where $\mu_{S^1}<\lambda$ of the lines
$x=\mathrm{constant}$, contributing $z=z_\delta$ to the superpotential. 
The two other families intersect the $y$-axis,
and can be described explicitly when $L_{r,\lambda}$ is a product torus 
(away from the blown up region): one consists as before of discs of radius
$r$ in the lines $y=\mathrm{constant}$, while the other one consists of
the proper transforms of discs which hit the $x$-axis at $(1,0)$, namely
the family of discs $z\mapsto (rz,\rho(rz-1)/(r-z))$ for fixed $|\rho|$.
Denote by $v$ the complexification of the right-pointing affine coordinate on 
Figure \ref{fig:blowup} in the chamber $r>1$, normalized so that,
if we ignore instanton corrections, the gluing across the wall 
$(r=1,\lambda>\epsilon)$ is given by $u=v^{-1}$. Then the two families
of discs intersecting the $y$-axis contribute respectively $v^{-1}$ and
$e^\epsilon z v^{-1}$ to the superpotential; the first family survives
the wall-crossing at $r=1$, while the second one degenerates by bubbling
of an exceptional disc (the part of the proper transform of the line $x=1$
where $\mu_{S^1}<\lambda$). This phenomenon is pictured on Figure
\ref{fig:blowup} (where the various discs are abusively represented as
tropical curves, which actually should be drawn in the complex affine
structure).

Thus the instanton-corrected gluing is
given by $u=v^{-1}+e^{\epsilon}zv^{-1}$ across the wall $(r=1,
\lambda>\epsilon)$; and a similar analysis shows that the portion of
the wall where $\lambda<\epsilon$ also gives rise to
the same instanton-corrected gluing. Thus, 
the instanton-corrected and completed mirror is given by
$$X^\vee=\{(u,v,z)\in \C^2\times\C^*,\ uv=1+e^{\epsilon}z\},\qquad W=u+z.$$
(Before completing the mirror by symplectically enlarging $X$, we would 
impose the restrictions $|u|<1$ and $|z|<1$.) The reader is referred to
\cite{AAK} for more details.
\medskip

\noindent {\bf Remark.} The above examples are particularly simple, as they
involve a single singularity of the special Lagrangian fibration and a single
wall-crossing correction. In more complicated examples, additional walls are
generated by intersections between the ``primary'' walls emanating from 
the singularities; in the end there are infinitely many walls, and hence
infinitely many instanton corrections to take into account when constructing
the mirror. A framework for dealing with
such situations has been introduced by Kontsevich and Soibelman \cite{KS2},
see also the work of Gross and Siebert \cite{GS1,GS2}.

\subsection{Beyond the Fano case: Hirzebruch surfaces}\label{ss:hirzebruch}

The construction of the mirror superpotential for toric Fano varieties
is well-understood (see e.g.\ \cite{hori,cho-oh,Au,FO3toric} for 
geometric derivations), and has been briefly summarized in \S \ref{ss:toric}
above. As pointed out to the author by Kenji Fukaya,
in the non-Fano case the superpotential differs from
the formula in Proposition \ref{prop:toric} by the presence of additional
terms, which count the virtual contributions of Maslov index 2 configurations
consisting of a disc of Maslov index 2 or more together with a collection
of spheres of non-positive Chern number. A non-explicit formula 
describing the general shape of the additional terms has been given
by Fukaya-Oh-Ohta-Ono: compare Theorems 3.4 and 3.5 in \cite{FO3toric}.
In this section we derive an {\it explicit} formula for the full
superpotential in the simplest example, using wall-crossing calculations.

The simplest non-Fano toric examples are rational ruled surfaces,
namely the Hirzebruch surfaces
$\mathbb{F}_n=\mathbb{P}(\mathcal{O}_{\PP^1}\oplus \mathcal{O}_{\PP^1}(n))$
for $n\ge 2$. The mirror of $\mathbb{F}_n$ is still $(\C^*)^2$, but with 
a superpotential of the form $W=W_0+\,$additional terms \cite{FO3toric},
where $W_0$ is given by (\ref{eq:toricW}), namely in this case
\begin{equation}\label{eq:W0}
W_0(x,y)=x+y+\frac{e^{-[\omega]\cdot[S_n]}}{xy^n}+
\frac{e^{-[\omega]\cdot[F]}}{y},
\end{equation}
where $[F]\in H_2(\mathbb{F}_n)$ is class of the fiber, and $[S_n]$ is the
class of a section of square $n$. 
The superpotential $W_0$ has $n+2$ critical points, four of which lie
within the region of $(\C^*)^2$ which maps to the moment polytope via
the logarithm map. Discarding the other critical points (i.e., restricting
to the appropriate subset of $(\C^*)^2$), homological mirror symmetry can 
be shown to hold for (a deformation of) $W_0$ \cite{AKO1} 
(see also \cite{abouzaid}). However, this is unsatisfactory for various
reasons, among others the discrepancy between the critical values of $W_0$
and the eigenvalues of quantum cup-product with the first Chern class
in $QH^*(\mathbb{F}_n)$ (see e.g.\ \S 6 of \cite{Au}, and \cite{FO3toric}).

\begin{figure}[t]
\setlength{\unitlength}{7mm}
\begin{picture}(5.3,2.5)(0,-0.25)
\psset{unit=\unitlength}
\psline(0,0)(3,0)(3,2)(0,2)(0,0)
\put(1.6,0.9){\small $\mathbb{F}_0$}
\put(1.4,-0.35){\tiny $0$}
\put(1.4,2.12){\tiny $0$}
\put(-0.28,0.85){\tiny $0$}
\put(3.1,0.85){\tiny $0$}
\psline{->}(3.8,1)(5,1)
\psline{->}(0.8,1)(0.3,1) \psline{->}(0.8,1)(0.8,0.5)
\put(0.4,1.15){\tiny $z_1$}
\put(0.9,0.8){\tiny $z_2$}
\end{picture}
\begin{picture}(5.5,2.5)(0,-0.25)
\psset{unit=\unitlength}
\pscurve(0,2)(0.2,0.6)(0.6,0.2)(3,0)
\pscurve[linestyle=dotted](0,2)(1.1,0.8)(3,0)
\psline(0,2)(3,2)(3,0)
\put(1.5,1.1){\small $\mathbb{F}_0$}
\put(1.4,2.12){\tiny $0$}
\put(3.1,0.85){\tiny $0$}
\put(1.1,0.8){\makebox(0,0)[cc]{\tiny $\times$}}
\put(-0.1,-0.1){\tiny $+2$}
\put(1.8,0.3){\tiny wall}
\psline{->}(0.7,1)(0.2,1) \psline{->}(0.7,1)(0.3,1.4)
\put(0.35,0.72){\tiny $u$}
\put(0.65,1.2){\tiny $w$}
\psline(3.8,1.05)(5,1.05)
\psline(3.8,0.95)(5,0.95)
\end{picture}
\begin{picture}(5.3,2.5)(0,-0.25)
\psset{unit=\unitlength}
\psline(0,0)(4,0)(1.2,2)
\pscurve(0,0)(0,0.5)(0,1)(0.25,1.8)(0.8,2)(1.2,2)
\pscurve[linestyle=dotted](0,0)(0.4,1.5)(1,1.6)(4,0)
\put(1.3,0.6){\small $\mathbb{F}_0$}
\put(-0.3,0.85){\tiny $0$}
\put(2.55,1.15){\tiny $0$}
\put(0.55,1.6){\makebox(0,0)[cc]{\tiny $\times$}}
\put(1.6,-0.33){\tiny $+2$}
\psline{->}(0.7,1)(0.7,0.5) \psline{->}(0.7,1)(0.3,0.6)
\put(0.8,0.72){\tiny $u$}
\put(0.2,0.88){\tiny $w$}
\psline{->}(3.8,1)(5,1)
\end{picture}
\begin{picture}(4,2.5)(0,-0.25)
\psset{unit=\unitlength}
\psline(0,2)(0,0)(4,0)(1,2)(0,2)
\put(1.3,0.6){\small $\mathbb{F}_2$}
\put(-0.3,0.85){\tiny $0$}
\put(2.55,1.15){\tiny $0$}
\put(0.2,2.1){\tiny $-2$}
\put(1.6,-0.33){\tiny $+2$}
\psline{->}(0.8,1)(0.3,1) \psline{->}(0.8,1)(0.8,0.5)
\put(0.5,1.15){\tiny $x$}
\put(0.9,0.8){\tiny $y$}
\end{picture}
\caption{Deforming $\mathbb{F}_0$ to $\mathbb{F}_2$}
\label{fig:f0f2}
\end{figure}
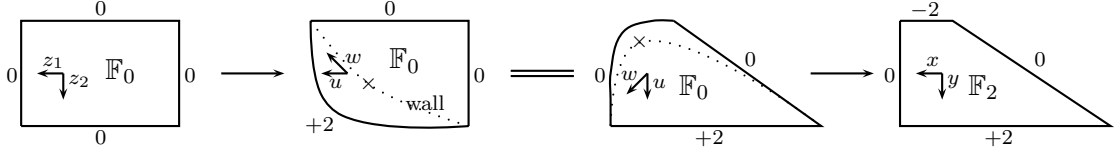

The approach we use to compute the full superpotential relies on the
observation that, depending on the parity of $n$, $\mathbb{F}_n$ is
deformation equivalent, and in fact symplectomorphic, to either 
$\mathbb{F}_0=\CP^1\times\CP^1$ or $\mathbb{F}_1$ (the one-point
blowup of $\CP^2$), equipped with a suitable symplectic form.
Carrying out the deformation explicitly provides a way of achieving
transversality for the Floer theory of Lagrangian tori in $\mathbb{F}_n$,
by deforming the non-regular complex structure of $\mathbb{F}_n$ to a
regular one. In the case of $\mathbb{F}_2$ and $\mathbb{F}_3$ at least,
the result of the deformation can be explicitly matched with $\mathbb{F}_0$
or $\mathbb{F}_1$ equipped with a non-toric holomorphic volume form of
the type considered in Example \ref{ex:conic}, which allows us to compute
the superpotential $W$. In the case of $\mathbb{F}_2$ the
deformation we want to carry out is pictured schematically in Figure \ref{fig:f0f2}.

\begin{proposition}\label{prop:f2}
The corrected superpotential on the mirror of\, $\F_2$ is the Laurent polynomial
\begin{equation}\label{eq:WF2}
W(x,y)=x+y+\frac{e^{-[\omega]\cdot[S_{+2}]}}{xy^2}+
\frac{e^{-[\omega]\cdot[F]}}{y}+
\frac{e^{-[\omega]\cdot[S_{-2}]}\,e^{-[\omega]\cdot[F]}}{y}.
\end{equation}
\end{proposition}

\noindent
This formula differs from $W_0$ by the addition of the last term;
geometrically, this term corresponds to configurations
consisting of a Maslov index 2 disc intersecting the exceptional section
$S_{-2}$ together with the exceptional section itself.

\proof
Consider the family of quadric surfaces
$\mathcal{X}=\{x_0x_1=x_2^2-t^2x_3^2\}\subset \CP^3\times\C$.
For $t\neq 0$, $X_t=\{x_0x_1=(x_2+tx_3)(x_2-tx_3)\}\subset \CP^3$ 
is a smooth quadric, and can be explicitly identified with the
image of the embedding of $\CP^1\times\CP^1$ given in homogeneous
coordinates by $$i_t:((\xi_0:\xi_1),(\eta_0:\eta_1))
\mapsto (\xi_0\eta_1:\xi_1\eta_0:\tfrac12(\xi_0\eta_0+\xi_1\eta_1):
\tfrac1{2t}(\xi_0\eta_0-\xi_1\eta_1)),$$
or, in terms of the affine coordinates $x=\xi_0/\xi_1$ and $y=\eta_0/\eta_1$,
$$i_t(x,y)=(x:y:\tfrac12(xy+1):\tfrac1{2t}(xy-1)).$$
For $t=0$, the surface $X_0=\{x_0x_1=x_2^2\}$ is a cone with vertex at the point
$(0\!:\!0\!:\!0\!:\!1)$, where the 3-fold $\mathcal{X}$ itself presents an ordinary
double point singularity. Denote by $\pi:\mathcal{X}'\to
\mathcal{X}$ a small resolution: composing with
the projection to $\C$, we obtain a family of surfaces $X'_t$, such that
$X'_t\cong X_t$ for $t\neq 0$, while $X'_0$ is the blowup of $X_0$, namely
$X'_0\cong \F_2$.

Consider the family of anticanonical divisors $$D_t=\{x_0=0,\ x_2=tx_3\}\cup
\{x_1=0,\ x_2=tx_3\}\cup \{x_3=-tx_2\}\subset X_t,$$ and equip $X_t$ with
a holomorphic volume form $\Omega_t$ with poles along $D_t$.
Observe that $D_t$ is the image by
$i_t$ of the lines at infinity $\xi_1=0$ and
$\eta_1=0$, and of the conic $C_t:\xi_0\eta_0=\frac{1-t^2}{1+t^2}\xi_1\eta_1$,
i.e.\ in affine coordinates, $xy=\frac{1-t^2}{1+t^2}$. Thus, for $t=1$
the pair $(X_t,D_t)$ corresponds to the toric anticanonical divisor in
$\CP^1\times\CP^1$, while for general $t$ the geometry of
$(X_t,D_t,\Omega_t)$ resembles closely that of Example \ref{ex:conic}
(the only difference is that we have compactified $\C^2$ to $\CP^1\times
\CP^1$). Finally, $D_0\subset X_0$ is precisely the toric anticanonical
divisor, consisting of the lines $x_0=0$, $x_1=0$ (two rays of the cone) 
and the conic $x_3=0$ (the base of the cone).

The quadrics $X_t$, the divisors $D_t$ and the volume forms $\Omega_t$
are preserved by the $S^1$-action $(x_0\!:\!x_1\!:\!x_2\!:\!x_3)\mapsto
(x_0e^{i\theta}\!:\!x_1e^{-i\theta}\!:\!x_2\!:\!x_3)$, and so is
the K\"ahler form induced by restriction of the Fubini-Study K\"ahler
form on $\CP^3$. Moreover, for $t=1$ the standard $T^2$-action on
$\CP^1\times\CP^1$ is induced by a subgroup of $PU(4)$, so that the
K\"ahler form is toric, and the configuration at $t=0$ is also toric
(with respect to a different $T^2$-action!); however for general $t$
we only have $S^1$-invariance.

The $S^1$-action lifts to the small resolution, and the lifted divisors
$D'_t=\pi^{-1}(D_t)\subset X'_t$ and holomorphic volume forms $\Omega'_t=
\pi^*\Omega_t$ are $S^1$-invariant. Additionally, $\mathcal{X}'$ can be
equipped with a $S^1$-invariant K\"ahler form, whose cohomology class 
depends on the choice of a parameter (the symplectic area of the 
exceptional $-2$-curve); restricting to $X'_t$, we obtain a family of
$S^1$-invariant K\"ahler forms $\omega'_t$. Moreover, careful choices
can be made in the construction in order to ensure that the K\"ahler
forms $\omega'_0$ and $\omega'_1$ on $X'_0$ and $X'_1$ are invariant 
under the respective $T^2$-actions.

  From the symplectic point of view, the family $(X'_t,\omega'_t)$ is
trivial; however, the complex structure for $t=0$ is non-generic.
At $t=1$ and $t=0$ the anticanonical divisors are the toric ones
for $\F_0$ and $\F_2$ respectively; deforming away from these values,
we partially smooth the toric anticanonical divisor, which in the
case of $\F_2$ requires a simultaneous deformation of the complex
structure because the exceptional curve $S_{-2}$ is rigid.
The deformation from $t=1$ to $t=0$ is now as pictured on
Figure~\ref{fig:f0f2}, which the reader is encouraged to keep in mind
for the rest of the argument.

For $t=1$, the mirror superpotential is given by the formula for the toric
Fano case (\ref{eq:toricW}), namely
\begin{equation}\label{eq:toricF0}
W=z_1+z_2+e^{-A}z_1^{-1}+e^{-B}z_2^{-1},\end{equation}
where $A$ and $B$ are the 
symplectic area of the two $\CP^1$ factors; the first two terms correspond
to discs contained in the affine chart with coordinates $x$ and $y$ we
have considered above, while the last two terms correspond to discs which
hit the lines at infinity. Deforming to general $t$, we have $X'_t\setminus
D'_t\cong \C^2\setminus \{xy=\frac{1-t^2}{1+t^2}\}$. Even though the
K\"ahler form $\omega'_t$ is not toric, the construction of an
$S^1$-invariant special Lagrangian fibration proceeds exactly as in Example
\ref{ex:conic}. The discussion carries over with only one modification:
when considered as submanifolds of $\CP^1\times\CP^1$,
the tori $T_{r,\lambda}=\{|xy-\frac{1-t^2}{1+t^2}|=r,\ \mu_{S^1}=\lambda\}$ 
bound additional families of Maslov index 2 holomorphic discs 
(intersecting the lines at infinity). 

In the chamber
$r>|\frac{1-t^2}{1+t^2}|$, deforming $T_{r,\lambda}$ to a product
torus shows that it bounds four families of Maslov index 2
discs, and the superpotential is given by (\ref{eq:toricF0}) as in
the toric case. On the other hand, in the chamber $r<|\frac{1-t^2}{1+t^2}|$, 
deforming $T_{r,\lambda}$ to the Chekanov torus shows that it bounds 
five families of Maslov index 2 holomorphic discs; explicit calculations
are given in Section 5.4 of \cite{Au}. (In \cite{Au} it was assumed for
simplicity that the two $\CP^1$ factors had equal symplectic areas, but
it is easy to check that the discussion carries over to the general case
without modification.) Using the same notations as
in Example \ref{ex:conic}, the superpotential is now given by 
\begin{equation}\label{eq:chekF0}
W=u+\frac{e^{-A}(1+w)}{uw}+\frac{e^{-B}(1+w)}{u}
\end{equation}
(see Corollary 5.13 in \cite{Au}). 
The first term $u$ corresponds to the family
of discs which are sections of $f:(x,y)\mapsto xy$ over the disc $\Delta$ of
radius $r$ centered at $\frac{1-t^2}{1+t^2}$; these discs pass through
the conic $xy=(1-t^2)/(1+t^2)$ and avoid all the toric divisors. 
The other terms correspond to sections of $f$ over $\CP^1\setminus \Delta$.
These discs intersect exactly one of the two lines at infinity, and one of
the two coordinate axes; each of the four possibilities gives rise to one
family of holomorphic discs. The various cases are as follows (see
Proposition 5.12 in \cite{Au}):

\begin{center}
\begin{tabular}{|c||c|c|c|c||c|}
\hline
class&$x=0$&$y=0$&$x=\infty$&$y=\infty$&weight\\
\hline
$H_1-\beta_0-\alpha$&no&yes&yes&no&$e^{-A\vphantom{A^2}}/uw$\\
$H_1-\beta_0$&yes&no&yes&no&$e^{-A}/u$\\
$H_2-\beta_0$&no&yes&no&yes&$e^{-B}/u$\\
$H_2-\beta_0+\alpha$&yes&no&no&yes&$e^{-B}w/u$\\
\hline
\end{tabular}
\end{center}
\noindent
Here $H_1=[\CP^1\times\{pt\}]$, $H_2=[\{pt\}\times\CP^1]$, and
$\beta_0$ and $\alpha$ are the classes in $\pi_2(\C^2,T_{r,\lambda})$
introduced in Example \ref{ex:conic}.

Here as in Example \ref{ex:conic},
it is easy to check that the two formulas for the superpotential
are related by the change of variables $w=z_1/z_2$ and $u=z_1+z_2$, which
gives the instanton-corrected gluing between the two chambers.

The tori $T_{r,\lambda}$ with $r<|\frac{1-t^2}{1+t^2}|$ cover the
portion of $X'_t\setminus D'_t$ where
$|xy-\frac{1-t^2}{1+t^2}|<|\frac{1-t^2}{1+t^2}|$, which under the
embedding $i_t$ corresponds to the inequality $$|x_2-tx_3|>
\left|\frac{2t}{1-t^2}\right|\,|tx_2+x_3|.$$
For $t\to 0$ this region covers almost all of $X'_t\setminus D'_t$, with
the exception of a small neighborhood of the lines $\{x_0=0,\ x_2=tx_3\}$
and $\{x_1=0,\ x_2=tx_3\}$. On the other hand, as $t\to 0$ the family of
special Lagrangian tori $T_{r,\lambda}$ converge to the standard toric
Lagrangian fibration on $X'_0=\F_2$, without any further wall-crossing 
as $t$ approaches zero provided that $r$ is small enough for $T_{r,\lambda}$
to lie within the correct chamber. It follows that, in suitable coordinates,
the superpotential for the Landau-Ginzburg mirror to $\F_2$ is given 
by (\ref{eq:chekF0}). 

All that remains to be done is to express the
coordinates $x$ and $y$ in (\ref{eq:WF2}) in terms
of $u$ and $w$. In order to do this, we investigate the limiting behaviors
of the five families of discs contributing to (\ref{eq:chekF0}) as 
$t\to 0$: four of these families are expected to converge to the
``standard'' families of Maslov index 2 discs in $\F_2$, since those are all
regular. Matching the families of discs allows us to match four of the
terms in (\ref{eq:chekF0}) with the four terms in $W_0$. The leftover
term in (\ref{eq:chekF0}) will then correspond to the additional term
in (\ref{eq:WF2}).

Consider a family of tori $T_{r(t),\lambda(t)}$
in $X'_t$ which converge to a $T^2$-orbit in $X'_0=\F_2$, corresponding
to fixed ratios $|x_0|/|x_3|=\rho_0$ and $|x_1|/|x_3|=\rho_1$ (and hence
$|x_2|/|x_3|=\sqrt{\rho_0\rho_1}$). Since the small resolution
$\pi:\mathcal{X}'\to \mathcal{X}$ is an isomorphism away from the
exceptional curve in $X'_0$, we can just work on $\mathcal{X}$ and use
the embeddings $i_t$ to convert back and forth between coordinates on
$X'_t\cong X_t\cong\CP^1\times\CP^1$ and homogeneous coordinates in $\CP^3$
for $t\neq 0$.
Since $$\frac{|x_2|}{|x_3|}=\left|t\,\frac{xy+1}{xy-1}\right|=\left|t+
\frac{2t}{xy-1}\right|$$
should converge to a finite non-zero value as $t\to 0$, the value of 
$|xy-1|$ must converge to zero, and hence
$r(t)=|xy-\frac{1-t^2}{1+t^2}|$ must also converge to $0$; in fact, an
easily calculation shows that
$r(t)\sim 2|t|/\sqrt{\rho_0\rho_1}$. Therefore, for $t$ small,
$xy$ is close to 1 everywhere on $T_{r(t),\lambda(t)}$. On the other hand,
$|x|/|y|=|x_0|/|x_1|$ converges to the finite value $\rho_0/\rho_1$;
thus $|x|$ and $|y|$ are bounded above and below on $T_{r(t),\lambda(t)}$.
Now, consider a holomorphic disc with boundary in $T_{r(t),\lambda(t)}$,
representing the class $H_1-\beta_0$. Since the $y$ coordinate has 
neither zeroes nor poles, by the maximum principle its norm is bounded
above and below by fixed constants (independently of $t$). The
point where the disc hits the line $x=\infty$ (i.e., $\xi_1=0$) has
coordinates $(x_0:x_1:x_2:x_3)=
(\xi_0\eta_1:0:\frac12 \xi_0\eta_0:\frac1{2t}\xi_0\eta_0)=(1:0:y:y/2t)$,
which given the bounds on $|y|$ converges to the singular point
$(0:0:0:1)$ as $t\to 0$. Thus, as $t\to 0$, this family of discs converges
to stable maps in $X'_0$ which have non-empty intersection with the 
exceptional curve. The same argument (exchanging $x$ and $y$) also applies
to the discs in the class $H_2-\beta_0$. On the other hand, the three other
families of discs can be shown to stay away from the exceptional curve.

In $\F_2$, the $T^2$-orbits bound four regular families of Maslov index 2
holomorphic discs, one for each component of the toric anticanonical divisor
$D'_0=S_{+2}\cup S_{-2}\cup F_0\cup F_1$; here $S_{-2}$ is the exceptional
curve, $S_{+2}$ is the preimage by $\pi$ of the component $\{x_3=0\}$ of
$D_0\subset X_0$,
and $F_0$ and $F_1$ are two fibers of the ruling, namely the proper
transforms under $\pi$ of the lines $x_0=0$ and $x_1=0$ in $X_0$.
The four families of discs can be constructed explicitly in coordinates
as in the proof of Proposition \ref{prop:toric}, see eq.\
(\ref{eq:holdisc}); regularity implies that, as we
deform $X'_0=\F_2$ to $X'_t$ for $t\neq 0$ small enough, all these discs
deform to holomorphic discs in $(X'_t,T_{r(t),\lambda(t)})$.

The term $y$ in $W_0$
corresponds to the family of discs intersecting the section $S_{+2}$, which
under the projection $X'_0\to X_0$ corresponds to the component $\{x_3=0\}$
of the divisor $D_0$). Thus, its deformation for $t\neq 0$ intersects
the component $\{x_3=-tx_2\}$ of the divisor $D_t$, namely the conic $C_t$
in the affine part of $\CP^1\times \CP^1$. Comparing the contributions to the
superpotential, we conclude that $y=u$. 

Next, the term $x$ in $W_0$ corresponds to the family of discs intersecting
the ruling fiber $F_0$, which projects to the line $\{x_0=0\}$ on $X_0$.
Thus, for small enough $t\neq 0$ these discs deform to a family of discs
in $X'_t$ that intersect the component $\{x_0=0,\ x_2=tx_3\}$ of $D'_t$,
i.e.\ the line at infinity $\eta_1=0$. There are two such families, in
the classes $H_2-\beta_0$ and $H_2-\beta_0+\alpha$; however we have seen
that the discs in the class $H_2-\beta_0$ approach the exceptional curve
as $t\to 0$, which would give a contradiction. Thus the term $x$ in $W_0$
corresponds to the family of discs in the class $H_2-\beta_0+\alpha$, which
gives $x=e^{-B}w/u$.

The proof is then completed by observing that the change of variables 
$x=e^{-B}w/u$, $y=u$ identifies (\ref{eq:WF2}) with (\ref{eq:chekF0}).
(Recall that the symplectic areas of $S_{+2}$, $S_{-2}$, and the ruling
fibers in $\F_2$ are respectively $A+B$, $A-B$, and $B$).

Note: as a quick consistency check, our change of variables matches the
term $e^{-(A+B)}/xy^2$ in (\ref{eq:WF2}), which corresponds to discs
in $\F_2$ that intersect the ruling fiber $F_1$, with the term $e^{-A}/uw$
in (\ref{eq:chekF0}), which corresponds to discs representing
the class $H_1-\beta_0-\alpha$ in $X'_t$ and intersecting the line
at infinity $\xi_1=0$. The remaining two terms in (\ref{eq:chekF0})
correspond to discs representing the classes $H_1-\beta_0$ and
$H_2-\beta_0$ in $X'_t$, whose limits as $t\to 0$ intersect the
exceptional curve $S_{-2}$, and can also be matched to the remaining
terms in (\ref{eq:WF2}). 
\endproof

A similar method can be applied to the case of $\F_3$, and yields:

\begin{proposition}\label{prop:f3}
The corrected superpotential on the mirror of\, $\F_3$ is the Laurent polynomial
\begin{equation}\label{eq:WF3}
W(x,y)=x+y+\frac{e^{-[\omega]\cdot[S_{+3}]}}{xy^3}+
\frac{e^{-[\omega]\cdot[F]}}{y}+
\frac{2e^{-[\omega]\cdot([S_{-3}]+2[F])}}{y^2}
+\frac{e^{-[\omega]\cdot([S_{-3}]+[F])}x}{y}.
\end{equation}
\end{proposition}

\begin{figure}[b]
\setlength{\unitlength}{7mm}
\begin{picture}(5,2.5)(0,-0.25)
\psset{unit=\unitlength}
\psline(0,0)(4,0)(2,2)(0,2)(0,0)
\put(1.7,0.8){\small $\mathbb{F}_1$}
\put(1.4,-0.35){\tiny $+1$}
\put(0.8,2.12){\tiny $-1$}
\put(-0.28,0.85){\tiny $0$}
\put(3.1,1.1){\tiny $0$}
\psline{->}(3.8,1)(4.8,1)
\psline{->}(0.8,1)(0.3,1) \psline{->}(0.8,1)(0.8,0.5)
\put(0.4,1.15){\tiny $z_1$}
\put(0.9,0.8){\tiny $z_2$}
\end{picture}
\begin{picture}(5.2,2.5)(0,-0.25)
\psset{unit=\unitlength}
\pscurve(0,2)(0.2,0.6)(0.8,0.2)(4,0)
\pscurve[linestyle=dotted](0,2)(1.1,0.8)(4,0)
\psline(0,2)(2,2)(4,0)
\put(1.5,1.1){\small $\mathbb{F}_1$}
\put(0.8,2.12){\tiny $-1$}
\put(3.1,1.1){\tiny $0$}
\put(1.1,0.8){\makebox(0,0)[cc]{\tiny $\times$}}
\put(0,-0.1){\tiny $+3$}
\put(1.8,0.3){\tiny wall}
\psline{->}(0.7,1)(0.2,1) \psline{->}(0.7,1)(0.3,1.4)
\put(0.35,0.72){\tiny $u$}
\put(0.65,1.2){\tiny $w$}
\psline(3.8,1.05)(4.8,1.05)
\psline(3.8,0.95)(4.8,0.95)
\end{picture}
\begin{picture}(5.5,2.5)(0,-0.25)
\psset{unit=\unitlength}
\psline(0,0)(5,0)(1.2,2)
\pscurve(0,0)(0,0.5)(0,1)(0.25,1.8)(0.8,2)(1.2,2)
\pscurve[linestyle=dotted](0,0)(0.4,1.5)(1,1.6)(5,0)
\put(1.8,0.6){\small $\mathbb{F}_1$}
\put(-0.3,1.85){\tiny $-1$}
\put(3.05,1.15){\tiny $0$}
\put(0.55,1.6){\makebox(0,0)[cc]{\tiny $\times$}}
\put(1.9,-0.33){\tiny $+3$}
\psline{->}(0.7,1)(0.7,0.5) \psline{->}(0.7,1)(0.3,0.6)
\put(0.8,0.72){\tiny $u$}
\put(0.2,0.88){\tiny $w$}
\psline{->}(4,1)(5.2,1)
\end{picture}
\begin{picture}(5,2.5)(0,-0.25)
\psset{unit=\unitlength}
\psline(0,2)(0,0)(5,0)(1,2)(0,2)
\put(1.8,0.6){\small $\mathbb{F}_3$}
\put(-0.3,0.85){\tiny $0$}
\put(3.05,1.15){\tiny $0$}
\put(0.2,2.1){\tiny $-3$}
\put(1.9,-0.33){\tiny $+3$}
\psline{->}(0.8,1)(0.3,1) \psline{->}(0.8,1)(0.8,0.5)
\put(0.5,1.15){\tiny $x$}
\put(0.9,0.8){\tiny $y$}
\end{picture}
\caption{Deforming $\mathbb{F}_1$ to $\mathbb{F}_3$}
\label{fig:f1f3}
\end{figure}
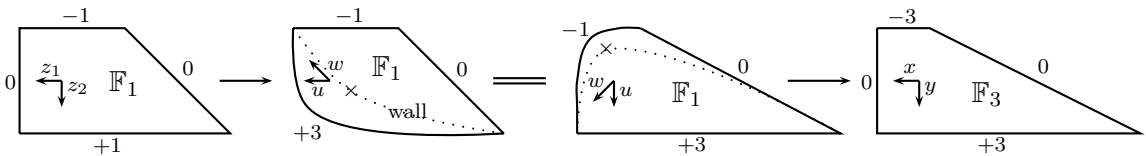

\proof[Sketch of proof]
The deformation we want to carry out is now depicted on Figure
\ref{fig:f1f3}.
One way of constructing this deformation is to start with the family
$\mathcal{X}'$ considered previously, and perform a birational
transformation. Namely, let $C'\subset\mathcal{X}'$ be the proper
transform of the curve $C=\{x_0=x_1=0,\ x_2=tx_3\}\subset\mathcal{X}$,
and let $\hat{\mathcal{X}}'$ be the blowup of $\mathcal{X}'$ along $C'$.
(This amounts to blowing up the point $x=y=\infty$ in each quadric $X'_t$
for $t\neq 0$, and the point where $S_{-2}$ intersects the fiber $F_1$ in
$X'_0\cong \F_2$). Next, let $\hat{Z}\subset \hat{\mathcal{X}}'$ be the
proper transform of the surface $Z=\{x_1=0,\ x_2=tx_3\}\subset\mathcal{X}$.
Denote by $\mathcal{X}''$ the 3-fold obtained by contracting $\hat{Z}$ in
$\hat{\mathcal{X}}'$: namely, $\mathcal{X}''$ is a family of surfaces
$X''_t$, each obtained from $X'_t$ by first blowing up a point as explained above
and then blowing down the proper transform of the line $x_1=0,\ x_2=tx_3$
(for $t\neq 0$ this is the line at infinity $\xi_1=0$, while for $t=0$ this
is the ruling fiber $F_1$). One easily checks that $X''_t\cong \F_1$ for
$t\neq 0$, while $X''_0\cong \F_3$.
Moreover, the divisors $D'_t\subset X'_t$ 
transform naturally under the birational transformations
described above, and yield a family of anticanonical divisors $D''_t\subset
X''_t$; for $t=0$ and $t=1$ these are precisely the toric anticanonical
divisors in $X''_0=\F_3$ and $X''_1=\F_1$.

In terms of the affine charts on $X'_t\cong \F_0$ considered in the proof of
Proposition \ref{prop:f2}, the birational transformations leading to
$(X''_t,D''_t)$ are performed ``at infinity'': thus $D''_t$ is again the
union of the conic $xy=\frac{1-t^2}{1+t^2}$ and the divisors at infinity,
namely for general $t$ we are again dealing with a compactified version
of Example \ref{ex:conic}. Thus $X''_t\setminus D''_t$ still contains an
$S^1$-invariant family of special Lagrangian tori $T_{r,\lambda}$, 
constructed as previously, and there are again two chambers separated by
the wall $r=|\frac{1-t^2}{1+t^2}|$; the only difference concerns the 
superpotential, since the tori $T_{r,\lambda}$ now bound
different families of holomorphic discs passing through the divisors at
infinity. These families and their contributions to the superpotential
can be determined by the same techniques as in the cases of $\CP^2$ and
$\CP^1\times\CP^1$, which are treated in Section 5 of \cite{Au}. Namely,
for $r>|\frac{1-t^2}{1+t^2}|$ the tori $T_{r,\lambda}$ can be isotoped
to product tori, and hence they bound four families of Maslov index 2 discs,
giving the familiar formula \begin{equation}
\label{eq:toricF1}
W=z_1+z_2+\frac{e^{-(A+B)}}{z_1z_2}+\frac{e^{-B}}{z_2},
\end{equation}
where $B$ is the area of the ruling fiber in $\F_1$ and $A$ is the area
of the exceptional curve. Meanwhile, for $r<|\frac{1-t^2}{1+t^2}|$ the
tori $T_{r,\lambda}$ can be isotoped to Chekanov tori; it can be shown
that they bound 6 families of Maslov index 2 discs, and using the same
notations as in Example \ref{ex:conic} we now have
\begin{equation}
\label{eq:chekF1}
W=u+\frac{e^{-(A+B)}(1+w)^2}{u^2w}+\frac{e^{-B}(1+w)}{u}.
\end{equation}
The instanton-corrected gluing between the two chambers is again 
given by $u=z_1+z_2$ and
$w=z_1/z_2$; in fact (\ref{eq:chekF1}) can be derived from
(\ref{eq:toricF1}) via this change of variables without having to
explicitly determine the holomorphic discs bounded by
$T_{r,\lambda}$. 

The strategy is now the same as in the proof of
Proposition \ref{prop:f2}: as $t\to 0$, the special Lagrangian 
fibrations on $X''_t\setminus D''_t$ converge to the standard 
fibration by $T^2$-orbits on $X''_0\simeq \F_3$, and the chamber 
$r<|\frac{1-t^2}{1+t^2}|$ covers arbitrarily large subsets of 
$X''_t\setminus D''_t$. Therefore, as before the superpotential for
the Landau-Ginzburg mirror to $\F_3$ is given by (\ref{eq:chekF1})
in suitable coordinates; the expression for the variables $x$ and $y$
in (\ref{eq:WF3}) in terms of $u$ and $w$ can be found by matching
some of the families of discs bounded by $T_{r,\lambda}$ as $t\to 0$
to the regular families of Maslov index 2 discs bounded by the $T^2$-orbits
in $\F_3$.

Concretely, the term $y$ in (\ref{eq:W0}) corresponds to holomorphic
discs in $\F_3$ which intersect the section $S_{+3}$. By regularity,
these discs survive the deformation to a small nonzero value of $t$,
and there they correspond to a family of discs which are entirely
contained in the affine charts. Hence, as before we must have $y=u$.
Identifying which term of (\ref{eq:chekF1}) corresponds to the term
$x$ in (\ref{eq:W0}) requires more work, but can be done exactly along
the same lines as for Proposition \ref{prop:f2}; in fact, we find
that it is given by exactly the same formula $x=e^{-B}w/u$ as in the
case of $\F_2$.
A posteriori this is not at all surprising, since this family of discs stays away
from the line at infinity $\xi_1=0$ in $X'_t$, and hence lies
in the part of $X'_t$ that is not affected by the birational transformations
that lead to $X''_t$.

Applying the change of variables $x=e^{-B}w/u$, $y=u$ to (\ref{eq:chekF1}),
and recalling that the symplectic areas of $S_{+3}$, $S_{-3}$ and the ruling
fibers in $\F_3$ are respectively $A+2B$, $A-B$, and $B$,
we arrive at (\ref{eq:WF3}), which completes the proof.
\endproof

It is tempting to interpret
the last two terms in (\ref{eq:WF3}) as the contributions of Maslov index 2
stable configurations that include the exceptional curve $S_{-3}$ as a
bubble component. Namely, the next-to-last term should be a virtual count
of configurations that consist of a double cover of a Maslov index 2
disc passing through $S_{-3}$, together with $S_{-3}$; and the last term
should be a virtual count of configurations consisting of a Maslov index 4
disc which intersects both the ruling fiber $F_0$ and the
exceptional section $S_{-3}$, together with $S_{-3}$.

In general, Fukaya-Oh-Ohta-Ono show that the ``naive'' superpotential
$W_0$ should be corrected by virtual contributions of Maslov index 2
configurations for which transversality fails in the toric setting;
moreover, they show that the perturbation data needed to make sense of
the virtual counts can be chosen in a $T^2$-equivariant manner
\cite{FO3toric}. In principle, different choices of perturbation data
could lead to different virtual counts of holomorphic discs, and hence
to different formulas for the corrected superpotential. Our approach
here can be understood as an {\it explicit} construction of a
perturbation that achieves transversality for holomorphic discs, by
deforming the complex structure to a generic one. However, our perturbation
is only $S^1$-equivariant rather than $T^2$-equivariant, so it is not
clear that our count of discs agrees with the virtual counts obtained
by using Fukaya-Oh-Ohta-Ono's perturbation data (the latter have not been
computed yet, in fact their direct computation seems extremely difficult). 
It is nonetheless our hope that the two counts might agree;
from this perspective it is encouraging to note that
open Gromov-Witten invariants are well-defined in the $S^1$-equivariant
setting, and not just in the toric setting \cite{liu}.

\subsection{Higher dimensions}

In this section we give two explicit local models for singularities of
Lagrangian fibrations in higher dimensions and their instanton-corrected
mirrors, generalizing the two examples considered in \S \ref{ss:examples}.
The open Calabi-Yau manifolds underlying the two examples are in fact
mirror to each other, as will be readily apparent. In complex dimension 3
these examples are instances of the two types of ``trivalent vertices'' 
that typically arise in the discriminant loci of special Lagrangian
fibrations on Calabi-Yau 3-folds and appear all over the relevant
literature (see e.g.\ \cite{grosstop}).  
These examples can also be understood by applying the
general machinery developed by Gross and Siebert \cite{GS1,GS2}; 
nonetheless, we find it interesting to have a fairly explicit and self-contained
description of the construction.

\subsubsection{} Consider $X=\C^n$, equipped with the standard K\"ahler form
$\omega$ and the holomorphic volume form $\Omega=(\prod x_i-\epsilon)^{-1}\,
dx_1\wedge\dots\wedge dx_n$, which has poles along the hypersurface
$D=\{\prod x_i=\epsilon\}$. Then $X\setminus D$ carries a fibration by
special Lagrangian tori $T_{r,\lambda}=\{(x_1,\dots,x_n)\in \C^n,\ 
|\prod x_i-\epsilon|=r,\ \mu_{T^{n-1}}(x_1,\dots,x_n)=\lambda\}$,
where $\mu_{T^{n-1}}:\C^n\to \R^{n-1}$ is the moment map for the
action of the group $T^{n-1}=\{\mathrm{diag}(e^{i\theta_1},\dots,
e^{i\theta_n}),\ \sum \theta_i=0\}$. More explicitly,
$$T_{r,\lambda}=\left\{(x_1,\dots,x_n)\in \C^n,\ \ 
\Bigl|{\textstyle\prod\limits_1^n x_i-\epsilon}\Bigr|=r,\ \ 
\tfrac12(|x_i|^2-|x_n|^2)=\lambda_i\ \forall i=1,\dots,n\!-\!1\right\}.$$
The tori $T_{r,\lambda}$ are $T^{n-1}$-invariant, and as in previous 
examples they are obtained by lifting special Lagrangian fibrations on the
reduced spaces.
As in Example \ref{ex:conic}, these tori are easiest to visualize
in terms of the projection $f:(x_1,\dots,x_n)\mapsto \prod x_i$, with
respect to which they fiber over circles centered at $\epsilon$; see
Figure \ref{fig:conic}. The main difference is that $f^{-1}(0)$
is now the union of the $n$ coordinate hyperplanes, and $T_{r,\lambda}$ is
singular whenever it hits the locus where the $T^{n-1}$-action
is not free, namely the points where at least two coordinates vanish.
Concretely, $T_{r,\lambda}$ is singular if and only if
$r=|\epsilon|$ and $\lambda$ lies in the {\it tropical
hyperplane} consisting of those $\lambda=(\lambda_1,\dots,\lambda_{n-1})$
such that either $\min(\lambda_i)=0$, or $\min(\lambda_i)$ is attained
twice. (For $n=3$ this is the union of the three half-lines $0=\lambda_1\le
\lambda_2$, $0=\lambda_2\le \lambda_1$, and $\lambda_1=\lambda_2\le 0$.)

By the maximum principle, any holomorphic disc in $(\C^n,T_{r,\lambda})$
which does not intersect $D=f^{-1}(\epsilon)$ must be contained inside
a fiber of $f$. The regular fibers of $f$ are diffeomorphic to
$(\C^*)^{n-1}$, inside which product tori do not bound any nonconstant
holomorphic discs. Hence, $T_{r,\lambda}$ bounds nontrivial
Maslov index 0 holomorphic discs if and only if $r=|\epsilon|$. In that
case, $T_{|\epsilon|,\lambda}$ intersects one of the components of
$f^{-1}(0)$ (i.e.\ a coordinate hyperplane isomorphic to $\C^{n-1}$)
in a product torus, which bounds various families of
holomorphic discs inside $f^{-1}(0)$.

The wall $r=|\epsilon|$ divides the moduli space of special Lagrangians 
into two chambers. In the chamber
$r>|\epsilon|$, the tori $T_{r,\lambda}$ can be be deformed into product
tori by a Hamiltonian isotopy that does not intersect $f^{-1}(0)$ (from
the perspective of the projection $f$, the isotopy amounts simply to deforming the
circle of radius $r$ centered at $\epsilon$ to a circle of the appropriate 
size centered at the origin). The product torus $S^1(r_1)\times\dots\times
S^1(r_n)$ bounds $n$ families of Maslov index 2 discs parallel
to the $x_1,\dots,x_n$ coordinate axes; denote their classes by 
$\beta_1,\dots,\beta_n$, and by
$z_i=\exp(-\int_{\beta_i}\omega)\,\hol_{\nabla}(\partial\beta_i)$ the
corresponding holomorphic weights. Thus we expect that $T_{r,\lambda}$
bounds $n$ families of Maslov index 2 holomorphic discs; these are
all sections of $f$ over the disc of radius $r$ centered at $\epsilon$,
and the discs in the class $\beta_i$ intersect the fiber $f^{-1}(0)$ at a
point of the coordinate hyperplane $x_i=0$. Since the deformation from
$T_{r,\lambda}$ to the product torus does not involve any wall-crossing,
the count of discs in the class $\beta_i$ is $1$, and the superpotential
is given by $W=z_1+\dots+z_n$.

Next we look at the chamber $r<|\epsilon|$. We first observe that the
Chekanov-type torus $T_{r,0}$ bounds only one family of Maslov index 2
holomorphic discs. Indeed, since Maslov index 2 discs have intersection
number 1 with $D=f^{-1}(\epsilon)$, they must be sections of $f$ over
the disc of radius $r$ centered at $\epsilon$, and hence they do not
intersect any of the coordinate hyperplanes. However, on $T_{r,0}$ we
have $|x_1|=\dots=|x_n|$, so the maximum principle applied to $x_i/x_n$
implies that the various coordinates $x_i$ are proportional to each other,
i.e.\ all such holomorphic discs must be contained in lines passing through
the origin. One easily checks that this gives a single family of holomorphic
discs; we denote by $\beta_0$ the corresponding homotopy class and by
$u=z_{\beta_0}$ the corresponding weight. Finally, since no exceptional
discs arise in the deformation 
of $T_{r,0}$ to $T_{r,\lambda}$,
we deduce that $T_{r,\lambda}$ also bounds a single
family of holomorphic discs in the class $\beta_0$, and that the
superpotential in the chamber $r<|\epsilon|$ is given by $W=u$.

When we increase the value of $r$ past $r=|\epsilon|$, with all
$\lambda_i>0$, the torus $T_{r,\lambda}$ crosses the coordinate
hyperplane $x_n=0$, and the family of holomorphic discs in the class
$\beta_0$ naturally deforms into the family of discs in the class $\beta_n$
mentioned above. However, the naive gluing $u=z_n$ must be corrected by
wall-crossing contributions. For $r=|\epsilon|$, $T_{r,\lambda}$ intersects
the hyperplane $x_n=0$ in a product torus. This torus bounds $n-1$ families
of discs parallel to the coordinate axes inside $\{x_n=0\}$, whose classes
we denote by $\alpha_1,\dots,\alpha_{n-1}$; we denote by $w_1,\dots,w_{n-1}$
the corresponding holomorphic weights, which satisfy $|w_i|=e^{-\lambda_i}$. It is easy to check that, on the
$r>|\epsilon|$ side, we have $\alpha_i=\beta_i-\beta_n$, and hence $w_i=z_i/z_n$;
general features of wall-crossing imply that $w_i$ should
not be affected by instanton corrections. Continuity of 
the superpotential across the wall implies that
the relation between $u$ and $z_n$ should be modified to
$u=z_1+\dots+z_n=z_n(w_1+\dots+w_{n-1}+1)$. Thus, only
the families of Maslov index 0 discs in the classes $\alpha_1,\dots,
\alpha_{n-1}$ contribute to the instanton corrections, even though
the product torus in $\{x_n=0\}$ also bounds higher-dimensional families
of holomorphic discs, whose classes are positive linear combinations of
the $\alpha_i$.

Similarly, when we increase the value of $r$ past $r=|\epsilon|$, with some
$\lambda_k=\min\{\lambda_i\}<0$, the torus  $T_{r,\lambda}$ crosses the
coordinate hyperplane $x_k=0$, and the family of discs in the class
$\beta_0$ deforms to the family of discs in the class $\beta_k$. However,
for $r=|\epsilon|$, $T_{r,\lambda}$ intersects the hyperplane $x_k=0$ in
a product torus, which bounds $n-1$ families of discs parallel to the
coordinate axes, representing the classes 
$\alpha_i-\alpha_k=\beta_i-\beta_k$ ($i\neq k,n$), with weight
$w_iw_k^{-1}=z_i/z_k$, and $-\alpha_k=\beta_n-\beta_k$, with weight
$w_k^{-1}=z_n/z_k$. The instanton-corrected gluing is now
$u=z_k(z_1/z_k+\dots+z_n/z_k+1)=z_1+\dots+z_n$.%

Piecing things together as in Example \ref{ex:conic}, we obtain a
description of
the corrected and completed SYZ mirror in terms of the
coordinates $u$, $v=z_n^{-1}$, $w_1,\dots,w_{n-1}$:

\begin{proposition} The mirror of $X=\C^n$ relatively to the divisor
$D=\{\prod x_i=\epsilon\}$ is
$$X^\vee=\{(u,v,w_1,\dots,w_{n-1})\in \C^2\times(\C^*)^{n-1},\ uv=1+w_1+
\dots+w_{n-1}\},\quad W=u.$$
\end{proposition}

A final remark: one way to check that the variables $w_i$ are indeed not affected
by the wall-crossing is to compactify $\C^n$ to $(\CP^1)^n$, equipped now
with the standard product K\"ahler form. Inside $(\CP^1)^n$ the tori 
$T_{r,\lambda}$ also bound families of Maslov index 2 discs that
pass through the divisors at infinity. These discs are sections of $f$
over the {\it complement} of the disc of radius $r$ centered at $\epsilon$,
and can be described explicitly in coordinates after deforming
$T_{r,\lambda}$ to either a product torus (for $r>|\epsilon|$) or a
Chekanov torus $T_{r,0}$ (for $r<|\epsilon|$). In the latter case, we
notice that the discs intersect the divisor at infinity once and $f^{-1}(0)$
once, so that in affine coordinates exactly one component of the map has a zero and
exactly one has a pole. Each of the $n^2$ possibilities gives one family
of holomorphic discs; the calculations are a straightforward adaptation of
the case of $\CP^1\times\CP^1$ treated in Section 5.4 of \cite{Au}.
The continuity of $W$ leads to an identity between the contributions 
to the superpotential coming from discs that intersect the compactification 
divisor ``$x_k=\infty$'' (a single family of discs for $r>|\epsilon|$, 
vs.\ $n$ families for $r<|\epsilon|$): namely, denoting by $\Lambda$ the
area of $\CP^1$, we must have
$$\frac{e^{-\Lambda}}{z_k}=\frac{e^{-\Lambda}}{uw_k}(w_1+\dots+w_{n-1}+1).$$
This is consistent with the formulas given above for the gluing between the
two chambers.

\subsubsection{} This example is treated carefully in \cite{AAK}, where
it is used as a standard building block to construct mirrors of
hypersurfaces in toric varieties. Here we only give an outline, 
for completeness
and for symmetry with the previous example.

Consider $\C^n$ equipped with the standard holomorphic volume form $\prod
d\log x_i$, and blow up the codimension 2 linear subspace $Y\times 0=\{
x_1+\dots+x_{n-1}=1,\ x_n=0\}$. This yields a complex manifold $X$ equipped with the
holomorphic volume form $\Omega=\pi^*(\prod d\log x_i)$, with poles along
the proper transform $D$ of the coordinate hyperplanes. The $S^1$-action rotating
the last coordinate $x_n$ lifts to $X$; consider an $S^1$-invariant K\"ahler
form $\omega$ for which the area of the $\CP^1$ fibers of the exceptional 
divisor is $\epsilon$ ($\epsilon\ll 1$), and which agrees with the standard K\"ahler form of
$\C^n$ away from a neighborhood of the exceptional divisor. 
Denote by $\mu_{S^1}:X\to \R$ the moment map of the
$S^1$-action, normalized to equal $0$ on the proper transform of the
coordinate hyperplane $x_n=0$, and $\epsilon$ at the stratum of fixed
points given by the section ``at infinity'' of the exceptional divisor.

The reduced spaces $X_\lambda=\{\mu_{S^1}=\lambda\}/S^1$ ($\lambda\ge 0$) are all 
smooth and diffeomorphic to $\C^{n-1}$. They carry natural holomorphic volume
forms, which are the pullbacks of $d\log x_1\wedge\dots\wedge d\log
x_{n-1}$, and K\"ahler forms $\omega_\lambda$. While $\omega_\lambda$
agrees with the standard K\"ahler form for $\lambda\gg \epsilon$, for
$\lambda<\epsilon$ the form $\omega_\lambda$ is not toric; rather, it
can be described as the result of collapsing a tubular neighborhood of size
$\epsilon-\lambda$ of the hypersurface $Y=\{x_1+\dots+x_{n-1}=1\}$ inside
the standard $\C^{n-1}$. Thus, it is not entirely clear that $X_\lambda$
carries a special Lagrangian torus fibration (though it does seem likely).

Nonetheless, using Moser's theorem to see that $\omega_\lambda$ is
symplectomorphic to the standard form on $\C^{n-1}$,
we can find a Lagrangian torus fibration on
the complement of the coordinate hyperplanes in 
$(X_\lambda,\omega_\lambda)$. Taking the preimages of these Lagrangians
in $\{\mu_{S^1}=\lambda\}$, we obtain a Lagrangian fibration on $X\setminus
D$, whose fibers are $S^1$-invariant Lagrangian tori $L_{r,\lambda}$;
for $\lambda\gg \epsilon$ these tori are of the form 
$$\{|\pi^*(x_i)|=r_i \ \,\forall 1\le i\le n-1,\ \ \ \mu_{S^1}=\lambda\}.$$

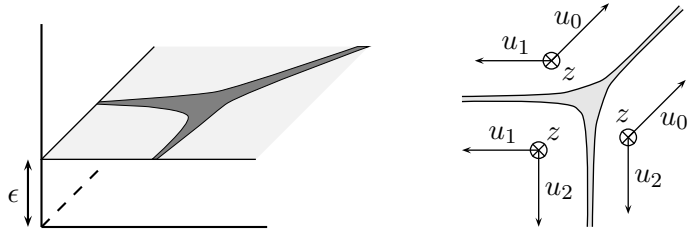
\begin{figure}[b]
\setlength{\unitlength}{1.5cm}
\begin{picture}(3.5,2)(0,0.2)
\psset{unit=\unitlength}
\newgray{vltgray}{0.95}
\newgray{ltgray}{0.9}
\newgray{dkgray}{0.85}
\pspolygon[linestyle=none,fillstyle=solid,fillcolor=vltgray](0,0.8)(1,1.8)(2.9,1.8)(1.9,0.8)
\psline[linestyle=dashed,dash=0.1 0.1](0,0.2)(0.5,0.7)
\pscustom[fillstyle=solid,fillcolor=gray,linewidth=0.3pt,curvature=.5 .1 0]{
  \psline(0.49,1.29)(0.9,1.27)
  \pscurve(0.9,1.27)(1.3,1.18)(0.98,0.8)
  \pscurve[liftpen=1](1.02,0.8)(1.4,1.1)(1.65,1.3)(2.9,1.8)
  \pscurve[liftpen=1](2.8,1.8)(1.55,1.4)(0.51,1.31)
}
\psline(0,0.2)(2,0.2)
\psline(0,0.2)(0,2)
\psline[linewidth=0.5pt](1,1.8)(0,0.8)(1.9,0.8)
\psline{<->}(-0.12,0.2)(-0.12,0.8) \put(-0.3,0.42){\small $\epsilon$}
\end{picture}
\setlength{\unitlength}{1.7cm}
\begin{picture}(1.75,1.8)(-1.2,-1)
\newgray{ltgray}{0.9}
\psset{unit=\unitlength,linewidth=0.5pt}
\pscustom[linestyle=none, fillstyle=solid,fillcolor=ltgray,curvature=.5 .1 0]{
  \pscurve[liftpen=1](0.7,0.73)(0.1,0.15)(-0.2,0.05)(-1,0.02)
  \pscurve[liftpen=1](-1,-0.02)(-0.12,-0.06)(-0.06,-0.12)(-0.02,-1)
  \pscurve[liftpen=1](0.02,-1)(0.05,-0.2)(0.15,0.1)(0.73,0.7)
}
\pscurve[curvature=.5 .1 0](0.7,0.73)(0.1,0.15)(-0.2,0.05)(-1,0.02)
\pscurve[curvature=.5 .1 0](-1,-0.02)(-0.12,-0.06)(-0.06,-0.12)(-0.02,-1)
\pscurve[curvature=.5 .1 0](0.02,-1)(0.05,-0.2)(0.15,0.1)(0.73,0.7)
\psline{->}(-0.4,-0.4)(-1,-0.4)
\psline{->}(-0.4,-0.4)(-0.4,-1)
\pscircle(-0.4,-0.4){0.07}
\psline(-0.35,-0.45)(-0.45,-0.35)
\psline(-0.35,-0.35)(-0.45,-0.45)
\put(-0.79,-0.31){\small $u_1$}
\put(-0.35,-0.75){\small $u_2$}
\put(-0.34,-0.32){\small $z$}
\psline{->}(-0.3,0.3)(-0.9,0.3)
\psline{->}(-0.3,0.3)(0.15,0.75)
\pscircle(-0.3,0.3){0.07}
\psline(-0.25,0.35)(-0.35,0.25)
\psline(-0.35,0.35)(-0.25,0.25)
\put(-0.69,0.39){\small $u_1$}
\put(-0.28,0.59){\small $u_0$}
\put(-0.22,0.15){\small $z$}
\psline{->}(0.3,-0.3)(0.3,-0.9)
\psline{->}(0.3,-0.3)(0.75,0.15)
\pscircle(0.3,-0.3){0.07}
\psline(0.25,-0.35)(0.35,-0.25)
\psline(0.35,-0.35)(0.25,-0.25)
\put(0.35,-0.65){\small $u_2$}
\put(0.55,-0.2){\small $u_0$}
\put(0.2,-0.18){\small $z$}
\end{picture}
\caption{$\C^3$ blown up along $\{x_1+x_2=1,\ x_3=0\}$}
\label{fig:blpants}
\end{figure}

The singularities of this fibration correspond to the fixed points of the
$S^1$-action inside $X\setminus D$, namely the ``section at infinity'' of
the exceptional divisor, defined by the equations
$\{\mu_{S^1}=\epsilon,\ \pi^*x_1+\dots+\pi^*x_{n-1}=1\}$. In the base
of the fibration, the discriminant locus is therefore of real codimension 1,
namely the amoeba of the hypersurface $Y$, sitting inside the affine
hyperplane $\lambda=\epsilon$ (see Figure \ref{fig:blpants} left). Moreover,
$L_{r,\lambda}$ bounds nonconstant discs of Maslov index 0 if and only
if it contains points where $\pi^*x_1+\dots+\pi^*x_{n-1}=1$. In that
case, the Maslov index 0 discs are contained in the total transforms of
lines parallel to the $x_n$-axis passing through a point of $Y\times 0$.
Thus, there are $n+1$ regions in which the tori $L_{r,\lambda}$ are weakly
unobstructed, corresponding to the connected components of the complement
of the amoeba of $Y$. 

To analyze holomorphic discs in $(X,L_{r,\lambda})$
and their contributions to the superpotential, we consider tori which
lie far away from the exceptional divisor and from the walls, i.e.\ for 
$r=(r_1,\dots,r_{n-1})$ sufficiently far from the amoeba of $Y$; then
$L_{r,\lambda}$ projects to a product torus in $\C^n$. When all $r_i\ll 1$
for all $i$, the maximum principle implies that holomorphic discs bounded
by $L_{r,\lambda}$ cannot hit the exceptional divisor; hence $L_{r,\lambda}$ 
bounds $n$ families of Maslov index 2 holomorphic discs, parallel to the 
coordinate axes. Denote by $\beta_1,\dots,\beta_{n-1},\delta$ the classes
of these discs, and by $u_1,\dots,u_{n-1},z$ the corresponding weights 
(i.e., the complexifications of the affine coordinates pictured in the
lower-left chamber of Figure \ref{fig:blpants} right).

Next consider the case where $r_k\gg 1$ and $r_k\gg r_i$ $\forall i\neq k$.
Then we claim that $L_{r,\lambda}$ now bounds $n+1$ families of Maslov index
2 holomorphic discs. Namely, since a Maslov index 2 disc intersects $D$ exactly once,
and the projections to the coordinates $(x_1,\dots,x_n)$ are holomorphic,
at most one of $\pi^*(x_1),\dots,\pi^*(x_{n-1})$ can be non-constant
over such a disc. Arguing as in the 2-dimensional case
(Example \ref{ex:c2blowup}), we deduce that $L_{r,\lambda}$ bounds $n$
families of discs parallel to the coordinate axes, and one additional
family, namely the proper transforms of Maslov index 4 discs in $\C^n$
which are parallel to the $(x_k,x_n)$-plane and hit the hyperplane $x_n=0$ at a point of $Y$.
Denote by $u_{1,(k)},\dots,u_{n-1,(k)},z_{(k)}$ the weights associated to 
the first $n$ families of discs: then the contribution of the
additional family to the superpotential is $e^{\epsilon}z_{(k)}u_{k,(k)}$.

Matching the contributions of the families of discs that intersect each
component of $D$, we conclude that the instanton-corrected gluings are given
by $z=z_{(k)}$, $u_i=u_{i,(k)}$ for $i\neq k$, and
$u_k=u_{k,(k)}(1+e^{\epsilon}z)$. Let $$u_{0,(k)}=\Bigl(\prod_{i=1}^{n-1}
u_{i,(k)}\Bigr)^{-1}=\Bigl(\prod_{i=1}^{n-1}u_i\Bigr)^{-1}(1+e^{\epsilon} z).$$
Then the coordinate $u_{0,(k)}$ is independent of $k$, and we can denote it
simply by $u_0$. The coordinates $(u_0,\dots,u_{n-1},z)$ can now be used to 
give a global description of the mirror (since forgetting one of the $u_i$
gives a set of coordinates for each chamber, as depicted in Figure
\ref{fig:blpants} right). Namely, after completion we arrive at:

\begin{proposition}[Abouzaid,---,Katzarkov \cite{AAK}]\label{prop:AAK}
The SYZ mirror of the blowup of\, $\C^n$ along
$\{x_1+\dots+x_{n-1}=1,\ x_n=0\}$ with anticanonical divisor the proper transform of
the toric divisor is
$$X^\vee=\{(u_0,\dots,u_{n-1},z)\in \C^n\times \C^*,\ u_0\dots u_{n-1}=
1+e^\epsilon z\},\quad W=u_1+\dots+u_{n-1}+z.$$
\end{proposition}

If instead we consider the blowup of $(\C^*)^{n-1}\times\C$ along the
{\it generalized pair of pants} $\{x_1+\dots+x_{n-1}=1,\ x_n=0\}$, i.e.\
we remove all the components of $D$ except the proper transform of the
$x_n=0$ coordinate hyperplane, then $X^\vee$ remains the same but the superpotential becomes simply
$W=z$ (since all the other terms in the above formula correspond to
discs that intersect the coordinate hyperplanes that we are now removing).

In \cite{AAK}, these local models are patched together in order to build
mirrors of more complicated blowups. The motivation for such a construction
comes from the observation that, if $Y$ is a hypersurface in $X$, then the derived category of $Y$ embeds into that
of the blowup of $X\times\C$ along $Y\times 0$ (this follows from a more general theorem
of Bondal and Orlov, see e.g.\ \cite{BO}); and, if $Y$
deforms in a pencil, then the Fukaya categories of these two
manifolds are also closely related (using Seidel's work; the key point is
that Lefschetz thimbles for a pencil in $X$ can be lifted to Lagrangian
spheres in the blowup of $X\times\C$ along $Y$). Thus, a mirror for the
blowup of $X\times\C$ along $Y$ is almost as good as a mirror for $Y$.
We illustrate this by considering one half of the homological mirror
symmetry conjecture in a very simple example.

Consider the case $n=3$ of Proposition \ref{prop:AAK} and its variants
where we remove various divisors from $D$. Consider the blowup of
$(\C^*)^2\times\C$ along $\{x_1+x_2=1,\ x_1,x_2\neq 0\}$ (a pair of pants,
i.e.\ $\PP^1$ minus three points): then $X^\vee$ is as in Proposition
\ref{prop:AAK}, i.e.\ (solving for $z$ as a function of $u_0,u_1,u_2$)
the complement of the hypersurface $u_0u_1u_2=1$ inside $\C^3$, and
the superpotential is $W=z=e^{-\epsilon}(u_0u_1u_2-1)$, whose critical
locus consists of the union of the three coordinate axes. Up to an
irrelevant scaling of the superpotential, this
Landau-Ginzburg model is indeed known to be a mirror to the pair of pants
(cf.\ work of Abouzaid and Seidel; see also \cite{Segenus2}). 
If instead we consider the blowup of
$\C^*\times\C^2$ along $\{x_1+x_2=1,\ x_1\neq 0\}$ ($\cong \C^*$), then
the superpotential becomes
$W=u_2+z=(e^{-\epsilon}u_0u_1+1)u_2-e^{-\epsilon}$; hence $W$ has a Morse-Bott
singularity along $M=\{u_0u_1=-e^\epsilon,\ u_2=0\}\simeq \C^*$, which
is mirror to $\C^*$. Finally, if we compactify
our example to consider the blowup of $\CP^2\times\C$ along a projective
line (given by $x_1+x_2=1$ in affine coordinates), then the mirror remains
the same manifold, but the superpotential acquires an extra term counting
discs that pass through the divisor at infinity, and becomes
$$W=e^{-\Lambda}u_0+u_1+u_2+z=e^{-\Lambda}u_0+u_1+u_2+e^{-\epsilon}u_0u_1u_2-e^{-\epsilon}$$
where $\Lambda$ is the area of a line in $\CP^2$.
This superpotential has two isolated non-degenerate critical points 
at $e^{-\Lambda}u_0=u_1=u_2=e^{\pm i\pi/2}e^{(\epsilon-\Lambda)/2}$, which is
reminiscent of the usual mirror of a $\CP^1$ with symplectic area
$\Lambda-\epsilon$ (to which our mirror can be
related by Kn\"orrer periodicity).

\section{Floer-theoretic considerations} \label{s:floer}

\subsection{Deformations and local systems}\label{ss:flavors}

There are at least three possible ways of deforming the Floer theory of
a given Lagrangian submanifold $L$ (for simplicity
we assume $L$ to be weakly unobstructed):

\begin{enumerate}
\item formally deforming the Floer theory of $L$ by an element $b\in CF^1(L,L)$;

\item equipping $L$ with a non-unitary local system;

\item deforming $L$ by a (non-Hamiltonian) Lagrangian isotopy and equipping
it with a unitary local system.
\end{enumerate}

Our goal in this paragraph is to explain informally how
these three flavors of deformation are related.
In particular, the careful reader will notice that Fukaya-Oh-Ohta-Ono
define the superpotential as a function on the moduli space of
{\it weak bounding cochains} for a given Lagrangian \cite{FO3book,FO3toric},
following the first approach, whereas in this paper and in \cite{Au} we view it
as a function on a moduli space of Lagrangians equipped with unitary local
systems, following the last approach.

Recall that there are several models for the Floer complex $CF^*(L,L)$.
We mostly consider the version in \cite{FO3book}, where the Floer
complex is generated by singular chains on $L$, representing incidence 
conditions at marked points on the boundary of holomorphic discs. The
$k$-fold product $\m_k$ is defined by 
\begin{equation}\label{eq:mk}
\m_k(C_1,\dots,C_k)=\!\!\!\sum_{\beta\in\pi_2(X,L)}\!\!\! z_\beta(L)\,
(ev_0)_*\left([\overline{\mathcal{M}}_{k+1}(L,\beta)]^{vir}\cap
ev_1^*C_1\cap\dots\cap ev_k^*C_k\right),\end{equation}
where $[\overline{\mathcal{M}}_{k+1}(L,\beta)]^{vir}$ is the (virtual)
fundamental chain of the moduli space of holomorphic discs in $(X,L)$
with $k+1$ boundary marked points representing the class $\beta$, 
$ev_0,\dots,ev_k$ are the evaluation maps at the marked points, and
$z_\beta$ is a weight factor as in (\ref{eq:zbeta}); when $k=1$ the
term with $\beta=0$ is replaced by a classical boundary term.

Here it is useful to also keep in mind a variant where
the Floer complex consists of differential forms or currents on $L$.
The product $\m_k$ is defined as in (\ref{eq:mk}), which now involves
pulling back
the given forms/currents to the moduli space of discs via the evaluation maps 
$ev_1,\dots,ev_k$ and pushing forward their product by integration
along the fibers of $ev_0$.
This setup allows us to ``smudge'' incidence conditions by replacing the 
integration current on a submanifold $C_i$
by a smooth differential form supported in a tubular neighborhood.

Given $b\in CF^1(L,L)$, Fukaya-Oh-Ohta-Ono \cite{FO3book} deform the 
$A_\infty$-algebra structure on the Floer complex by setting
\begin{equation}\label{eq:mkb}
\m_k^b(C_1,\dots,C_k)=\sum_{l=l_0+\dots+l_k\ge 0}
\m_{k+l}(\,\underbrace{b,\dots,b}\limits_{l_0}\,,C_1,\,
\underbrace{b,\dots,b}\limits_{l_1}\,,\dots,C_k,\,\underbrace{b,\dots,b}
\limits_{l_k}\,).\end{equation}
We will actually restrict our attention to the case where $b$ is a 
{\it cycle}, representing a class $[b]\in H^1(L)$ (or, dually, in $H_{n-1}(L)$).

Working over the Novikov ring, the sum (\ref{eq:mkb}) is guaranteed to be
well-defined when $b$ has coefficients in the maximal ideal
\begin{equation}\label{eq:novikov+}
\Lambda_+=\left\{\textstyle\sum\limits_i\, a_i\, T^{\lambda_i}\in\Lambda_0\,\Big|\,
\lambda_i>0\right\}\end{equation} of
$\Lambda_0=\{\textstyle\sum a_i\, T^{\lambda_i}\,|\,
a_i\in\C,\ \lambda_i\in\R_{\ge 0},\ \lambda_i\to +\infty\}$.
However, it has been observed by Cho~\cite{cho-nonunitary} (see also
\cite{FO3toric}) that, in the toric case, the sum (\ref{eq:mkb})
is convergent even when $b$ is a general element of $H^1(L,\Lambda_0)$.
Similarly, in favorable cases (at least for toric Fanos) we can 
also hope to make sense of (\ref{eq:mkb}) when working over $\C$ (in the
``convergent power series'' setting); however in general this poses
convergence problems.

The second type of deformation we consider equips $L$ with a local
system (a flat connection), characterized by its holonomy $\hol_\nabla$,
which is a homomorphism from $\pi_1(L)$ to $\Lambda_0^*$ (the multiplicative
group formed by elements of the Novikov ring with nonzero coefficient of
$T^0$)  or $\C^*$. The local system modifies the weight $z_\beta$
for the contribution to $\m_k$ of discs in the class $\beta$
by a factor of $\hol_\nabla(\partial\beta)$.

\begin{lemma}\label{l:deformequiv}
For any cycle $b$ such that convergence holds, the deformation 
of the $A_\infty$-algebra $CF^*(L,L)$ given by (\ref{eq:mkb}) is 
equivalent to equipping $L$ with a local system with holonomy $\exp(b)$,
i.e.\ such that
$\hol_\nabla(\gamma)=\exp([b]\cdot [\gamma])$ for all $\gamma\in \pi_1(L)$.
\end{lemma}

\proof[Sketch of proof]
The statement reduces to a calculation showing that, given a holomorphic
disc $u\in \mathcal{M}_{k+1}(L,\beta)$ (or more generally an element of
the compactified moduli space),
the contribution of ``refined'' versions of $u$ 
(with extra marked points mapped to $b$) to $\m_k^b$ is 
$\exp([b]\cdot[\partial\beta])$ times the contribution of $u$ to $\m_k$. 

This is easiest to see when we represent
the class $[b]$ by a smooth closed 1-form on $L$.
For fixed $\underline{l}=(l_0,\dots,l_k)$,
consider the forgetful map $\pi_{\underline{l}}: 
\mathcal{M}_{k+l+1}(L,\beta)\to \mathcal{M}_{k+1}(L,\beta)$
which deletes the marked points corresponding to the $b$'s in
(\ref{eq:mkb}), and its extension $\bar\pi_{\underline{l}}$ to the
compactified moduli spaces.
The fiber of $\pi_{\underline{l}}$ above 
$u\in\mathcal{M}_{k+1}(L,\beta)$ is
a product of open simplices of dimensions $l_0,\dots,l_k$, parametrizing
the positions of the $l_0+\dots+l_k$ new marked points along the intervals
separated by the $k+1$ marked points of $u$ on the boundary 
of the disc; we denote by
$\Delta_{\underline{l}}$ the corresponding subset of $(\partial D^2)^l$.
The formula for $\m_{k+l}(b^{\otimes l_0},C_1,b^{\otimes l_1},\dots,C_k,
b^{\otimes l_k})$ involves an integral over
$\overline{\mathcal{M}}_{k+l+1}(L,\beta)$,
but this integral can be pushed
forward to $\overline{\mathcal{M}}_{k+1}(L,\beta)$
by integrating over the fibers of $\bar{\pi}_{\underline{l}}$;
the resulting integral differs from that for $\m_k(C_1,\dots,C_k)$
by an extra factor
$\int_{\bar\pi^{-1}_{\underline{l}}(u)} \prod ev_i^*b =
\int_{\smash{\,\overline{\Delta}_{\underline{l}}}}  \prod
(u_{|\partial D^2}\circ pr_i)^*b$ in the integrand. 

Note that this calculation assumes that the
virtual fundamental chains have been constructed consistently, so that 
$[\overline{\mathcal M}_{k+l+1}(L,\beta)]^{vir}=\bar\pi_{\underline{l}}^*
([\overline{\mathcal M}_{k+1}(L,\beta)]^{vir})$ as expected.
Achieving this property is in general a non-trivial problem.

Next we sum over $\underline{l}$: the subsets
$\overline{\Delta}_{\underline{l}}$ of $(\partial D^2)^l$ have disjoint
interiors, and their union $\overline{\Delta}$ is the set of all $l$-tuples
of points which lie in counterclockwise order on the interval obtained
by removing the outgoing marked point of $u$ from $\partial D^2$.
By symmetry, the integral of $\prod (u_{|\partial D^2}\circ pr_i)^*b$
over $\overline{\Delta}$ is $1/l!$ times the integral over $(\partial
D^2)^l$. Thus
$$\sum_{\underline{l}} \int_{\overline{\Delta}_{\underline{l}}}
\,{\textstyle\prod\limits_{i=1}^l}\, (u\circ pr_i)^*b\,=\,
\frac{1}{l!}\int_{(\partial D^2)^l} 
\,{\textstyle\prod\limits_{i=1}^l}\, (u\circ pr_i)^*b
= \frac{1}{l!}\Bigl(\int_{\partial D^2}
u^*b\Bigr)^l=\frac{([b]\cdot[\partial\beta])^l}{l!}.$$
The statement then follows by summing over $l$.
\endproof

One can also try to prove Lemma \ref{l:deformequiv} working entirely
with chains on $L$ instead of differential forms, but it is technically
harder. If we take $b$ to be a codimension 1 cycle in $L$ and attempt to
reproduce the above argument, the incidence constraints at the additional
marked points (all mapping to $b$) are not transverse to each other. 
In fact, $\m_k^b$ will include contributions from stable maps with
constant disc bubbles mapping to $b$. The difficulty is then to 
understand the combinatorial rule for counting such contributions,
or more precisely, why a constant bubble with $j$
marked points on it, all mapped to a same point of $b$, should contribute
a combinatorial factor of $1/j!$. 

The equivalence between the two types of deformations also holds if we
consider not just $L$ itself, but the whole Fukaya category. Given
a collection of Lagrangian submanifolds $L_0,\dots,L_k$ with $L_{i_0}=L$
for some $i_0$, the Floer theoretic product
$\m_k:CF^*(L_0,L_1)\otimes\dots\otimes CF^*(L_{k-1},L_k)\to CF^*(L_0,L_k)$
can again
be deformed by a cycle $b\in CF^1(L,L)$. Where the usual product $\m_k$ is a
sum over holomorphic discs with $k+1$ marked points, the deformed product
$\m_k^b$ counts
discs with an arbitrary number of additional marked points, all lying on the
interval of $\partial D^2$ which gets mapped to $L$, and with
inputs $b$ inserted accordingly into the Floer product as in (\ref{eq:mkb}).
By the same argument as above, if we represent $b$ by a closed 1-form on
$L$, and consider discs with fixed corners and in a fixed homotopy class 
$\beta$, the deformation amounts to the insertion of an extra factor
$\exp(\int_{\partial \beta\cap L} b)$. Meanwhile, equipping $L$ with
a flat connection $\nabla$ affects the count of discs in the class
$\beta$ by a factor $\hol_\nabla(\partial\beta\cap L)$. Thus, if we
ensure that the two match, e.g.\ by choosing $\nabla=d+b$, the two
deformations are again equivalent.
\medskip

Next, we turn to the relation between non-unitary local systems and
non Hamiltonian deformations. Consider a deformation of $L$ to a nearby
Lagrangian submanifold $L_1$; identifying a tubular neighborhood of $L$
with a neighborhood of the zero section in $T^*L$, we can think of
$L_1$ as the graph of a $C^1$-small closed form $\varphi\in\Omega^1(L,\R)$.
Assume that $L$ can be isotoped to $L_1$ (e.g.\ through
$L_t=\mathrm{graph}(t\varphi)$) in such a way that there is a one-to-one
correspondence between the holomorphic discs bounded by $L$ and $L_1$,
namely given a class $\beta\in \pi_2(X,L)$ and the corresponding class
$\beta_1\in \pi_2(X,L_1)$, we have $\mathcal{M}_{k+1}(L,\beta)\simeq
\mathcal{M}_{k+1}(L_1,\beta_1)$. Observing that $\int_{\beta_1}\omega=
\int_\beta\omega+\int_{\partial\beta}\varphi$, deforming $L$ to $L_1$
affects the contribution of these discs by a factor of $\exp(-[\varphi]
\cdot[\partial\beta])$. Thus, deforming $L$ to $L_1$ is equivalent to 
equipping $L$ with a local system with holonomy $\exp(-[\varphi])$;
for example we can set $\nabla=d-\varphi$. (This is when working over
complex numbers; over the Novikov ring we would similarly want to
equip $L$ with a local system with holonomy $T^{[\varphi]}$.)
However, this only works as long as there is a good correspondence
between moduli spaces of holomorphic discs; while the assumption we
made can be relaxed to some extent, we cannot expect things to work
so simply when the deformation from $L$ to $L_1$ involves wall-crossing.

Similarly, given another Lagrangian submanifold $L'$, if the isotopy 
from $L$ to $L_1=\mathrm{graph}(\varphi)$ can be carried out in a
manner that remains transverse to $L'$ at all times then we can hope to define
a map from $CF^*((L,\nabla),L')$ (with $\hol(\nabla)=\exp(-[\varphi])$)
to $CF^*(L_1,L')$ in a manner compatible
with all Floer-theoretic products as long as we can find a one-to-one
correspondence between the relevant holomorphic discs. One could also
try to proceed in a slightly greater degree of generality by attempting
to construct continuation maps between the Floer complexes (as in
the usual proof of Hamiltonian isotopy invariance of Floer homology).
However, one should keep in mind that this is doomed to fail in general.
For instance, consider $X=S^2$, let $L_1$ be the equator, and $L$ a circle
parallel to the equator but disjoint from it. Denote by $A$ the 
annulus bounded by $L$ and $L_1$, and equip $L$ with a non-unitary
local system $\nabla$ with holonomy $\exp(\int_A\omega)$. One easily checks that
the Lagrangians $L_1$ and $(L,\nabla)$ have well-defined and non-vanishing Floer homology,
and the  $A_\infty$-algebras $CF^*(L_1,L_1)$ and 
$CF^*((L,\nabla),(L,\nabla))$ are isomorphic (by the
argument above). However, $CF^*((L,\nabla),L_1)=0$ since $L$ and $L_1$
are disjoint, so $(L,\nabla)$ and $L_1$ cannot be isomorphic. 
(See also the discussion in \S \ref{ss:failure}).
\medskip

\noindent {\bf Remark.} Specializing (\ref{eq:mkb}) to $k=0$,
the identity $\m_0^b=\m_0+\m_1(b)+\m_2(b,b)+\dots$
offers a simple perspective into the idea that the derivatives
of the superpotential $W$ at a point
$\mathcal{L}=(L,\nabla)$ encode information about the (symmetrized)
Floer products $\m_k$ on $CF^*(\mathcal{L},\mathcal{L})$,
as first shown by Cho in \cite{Cho}. In particular, one can re-derive
from this identity the fact that, if $\mathcal{L}$ is not a
critical point of the superpotential, then the fundamental class
of $L$ is a Floer coboundary and $HF^*(\mathcal{L},\mathcal{L})$
vanishes. (For a direct proof, see \cite{cho-oh,Cho},
see also \S 6 of \cite{Au}.)

\subsection{Failure of invariance and divergence issues}\label{ss:failure}

In this section, we look more carefully into a subtle issue with instanton
corrections and the interpretation of the mirror as a moduli space of
Lagrangian submanifolds up to Floer-theoretic equivalence. We return to
Example \ref{ex:conic}, i.e.\ $\C^2$ equipped with the standard K\"ahler
form and the holomorphic volume form $\Omega=dx\wedge dy/(xy-\epsilon)$,
and use the same notations as above. Consider two
special Lagrangian fibers on opposite sides of the wall, $T_1=T_{r_1,0}$ and 
$T_2=T_{r_2,0}$, where $r_1<|\epsilon|<r_2$ are chosen in a way such that
the points of $M$ corresponding to $T_1$ and $T_2$ (equipped with the
trivial local systems) are identified under the instanton-corrected gluing
$u=z_1+z_2$. Namely, the torus $T_1$ corresponds to a point with
coordinates  $w=1$, $u=\exp(-A_1)\in \R_+$, where $A_1$ is the symplectic area
of a Maslov index 2 disc in $(\C^2,T_1)$, e.g.\ either of the two portions
of the line $x=y$ where $|xy-\epsilon|<r_1$; meanwhile, $T_2$ corresponds to
$w=1$, $z_1=z_2=\exp(-A_2)\in \R_+$, where $A_2$ is the symplectic area of
a Maslov index 2 disc in $(\C^2,T_2)$, or equivalently half of the area
of the portion of the line $x=y$ where $|xy-\epsilon|<r_2$. 
The area $A_i$ can be expressed by an explicit formula in terms of
$r_i$ and $\epsilon$; the actual relation is irrelevant, all that matters
to us is that $A_i$ is a monotonically increasing function of $r_i$.
Now we choose 
$r_1$ and $r_2$ such that $\exp(-A_1)=2\exp(-A_2)$ and $r_1<|\epsilon|<r_2$.

We will consider the tori $T_1$ and $T_2$ inside $X^0=X\setminus
D=\C^2\setminus \{xy=\epsilon\}$, where they do not bound any nonconstant
holomorphic discs. (Another option would be to instead compactify
$\C^2$ to $\CP^2$, and choose the parameters of the construction so that
$\exp(-A_1)=2\exp(-A_2)=\exp(-\frac13 \int_{\CP^1}\omega)$; then $T_1$
and $T_2$ would be weakly unobstructed and would still have non-vanishing
convergent power series Floer homology. The discussion below would carry over
with minor modifications.)

Working in $X^0$, the convergent power series Floer homologies 
$HF^*(T_1,T_1)$ and $HF^*(T_2,T_2)$ are isomorphic to each other
(and to the cohomology of $T^2$). In fact the same property would hold
for any other $T_{r,\lambda}$ due to the absence of holomorphic discs in
$X^0$, but
in the case of $T_1$ and $T_2$ we expect to have a {\it distinguished}\/
isomorphism between
the Floer homology groups,
considering that $T_1$ and $T_2$ are in the same instanton-corrected
equivalence class and meant to be ``isomorphic''. 
However, $T_1$ and $T_2$ are disjoint, so $CF^*(T_1,T_2)$ is zero, which
does not allow for the existence of the expected isomorphism.
(Note that the issue would not arise when working over the Novikov ring:
we would then have needed to choose the areas $A_1$ and $A_2$ above so that
$T^{-A_1}=2T^{-A_2}$, which never holds. In that case, one should instead
take $A_1=A_2$ and equip $T_1$ with a nontrivial local system; but then
$T_1$ and $T_2$ cannot be made disjoint by Hamiltonian isotopies.)

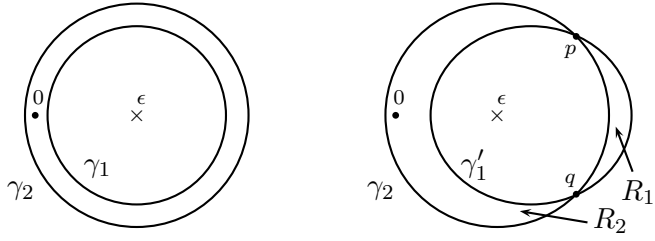
\begin{figure}[t]
\setlength{\unitlength}{1.5cm}
\begin{picture}(2,2)(-1,-1)
\psset{unit=\unitlength}
\pscircle(0,0){1}
\pscircle(0,0){0.8}
\pscircle*(-0.9,0){0.03}
\put(0,0){\makebox(0,0)[cc]{\tiny $\times$}}
\put(-0.92,0.1){\tiny $0$}
\put(0,0.1){\tiny $\epsilon$}
\put(-0.47,-0.5){\small $\gamma_1$}
\put(-1.15,-0.7){\small $\gamma_2$}
\end{picture}
\qquad\qquad
\begin{picture}(2,2)(-1,-1)
\psset{unit=\unitlength}
\pscircle(0,0){1}
\psellipse(0.3,0)(0.9,0.8)
\pscircle*(-0.9,0){0.03}
\put(0,0){\makebox(0,0)[cc]{\tiny $\times$}}
\put(-0.92,0.1){\tiny $0$}
\put(0,0.1){\tiny $\epsilon$}
\put(-0.32,-0.5){\small $\gamma'_1$}
\put(-1.15,-0.7){\small $\gamma_2$}
\psline{->}(0.8,-0.95)(0.2,-0.85) \put(0.85,-1){\small $R_2$}
\psline{->}(1.2,-0.5)(1.05,-0.1) \put(1.1,-0.75){\small $R_1$}
\pscircle*(0.7,-0.7){0.03} \put(0.6,-0.6){\tiny $q$}
\pscircle*(0.7,0.7){0.03} \put(0.6,0.55){\tiny $p$}
\end{picture}

\caption{Creating intersections between $T_1$ and $T_2$} \label{fig:wiggle}
\end{figure}

One way to realize the isomorphism between $T_1$ and $T_2$ is to deform
one of them by a Hamiltonian isotopy (without crossing any walls) in order
to create intersections. Namely, projecting $\C^2$ to $\C$ by the map
$f(x,y)=xy$, $T_1$ and $T_2$ fiber above concentric circles $\gamma_i=
\{|z-\epsilon|=r_i\}$, and inside each fiber they consist of the 
``equatorial'' $S^1$-orbit where $|x|=|y|$. Deform $T_1$
by a Hamiltonian isotopy, without crossing $\epsilon$ or $0$, to a
$S^1$-invariant Lagrangian torus $T'_1$ which fibers above a closed 
curve $\gamma'_1$ intersecting $\gamma_2$ in two points $p$ and $q$,
and $T'_1=f^{-1}(\gamma'_1)\cap\{|x|=|y|\}$ (see Figure \ref{fig:wiggle}).
Then $T'_1$ and $T_2$ intersect along two circles, which can be handled
either as a degenerate Morse-Bott type intersection ($CF^*(T'_1,T_2)$ is
then generated by chains on $T'_1\cap T_2$), or by further
perturbing $T'_1$ to replace each $S^1$ by two transverse
intersection points.

\begin{proposition}
In $X^0=\C^2\setminus f^{-1}(\epsilon)$, the convergent power series
Floer homology $HF^*(T'_1,T_2)$ is well-defined and isomorphic to $H^*(T^2,\C)$.
\end{proposition}

\proof
Any holomorphic disc in $X^0=\C^2\setminus f^{-1}(\epsilon)$ that contributes
to the Floer differential on $CF^*(T'_1,T_2)$ is necessarily a section of
$f$ over one of the two regions $R_1$ and $R_2$ delimited by $\gamma'_1$
and $\gamma_2$ (see Figure \ref{fig:wiggle}). Recalling that $|x|=|y|$ on
$T'_1\cup T_2$, the maximum principle applied to $x/y$ implies that,
if a disc with boundary in $T'_1\cup T_2$ intersects neither the $x$ axis
nor the $y$ axis, then $x/y$ is constant over it. Thus, there is 
exactly one $S^1$-family of such sections of $f$ over $R_1$, namely the
portions of the lines $y=e^{i\varphi}x$ which lie in $f^{-1}(R_1)$.
On the other
hand, there are two $S^1$-families of sections over $R_2$. Indeed, 
let $g:D^2\to R_2$ be a biholomorphism given by the Riemann mapping 
theorem, chosen so that $g(0)=0$, and consider
a holomorphic map $u:(D^2,\partial D^2)\to (X^0,T'_1\cup T_2)$,
$z\mapsto u(z)=(x(z),y(z))$ such
that $f\circ u$ maps $D^2$ biholomorphically onto $R_2$. Up to a reparametrization
we can assume that $f\circ u=g$. Over the image of $u$, either
$x$ or $y$ must vanish transversely once; assume that it is $x$ that
vanishes. Then $z\mapsto x(z)/y(z)$ is a holomorphic function on the disc,
taking values
in the unit circle along the boundary, and vanishing once at the origin,
therefore it is of the form $z\mapsto e^{i\varphi}z$ for some
$e^{i\varphi}\in S^1$. Thus
$u(z)=(e^{i\varphi/2}(zg(z))^{1/2},e^{-i\varphi/2}(g(z)/z)^{1/2})$.
This gives an $S^1$-family of holomorphic sections over $R_2$; the other
one is obtained similarly by exchanging $x$ and $y$.

Denote by $\alpha_1$ (resp.\ $\alpha_2$) the symplectic area of the 
holomorphic discs in $(X^0,T'_1\cup T_2)$ which are sections of $f$
over $R_1$ (resp.\ $R_2$). By construction, these areas are related to
those of the Maslov index 2 discs bounded by $T'_1$ and $T_2$ in $\C^2$:
namely, $\alpha_2-\alpha_1=A_2-A_1$. Thus, the choices made above imply
that $\exp(-\alpha_1)=2\exp(-\alpha_2)$. After a careful check of signs,
this in turn implies that the contributions of the various holomorphic 
discs in $(X^0,T'_1\cup T_2)$ to the Floer differential on 
$CF^*(T'_1,T_2)$ (with $\C$ coefficients) cancel out.
\endproof

\noindent
Denote by $e_p$ the generator of $CF^0(T'_1,T_2)$ which comes
from the intersections in $f^{-1}(p)$, and denote by $e_q$ the generator
of $CF^0(T_2,T'_1)$ which comes from the intersections in $f^{-1}(q)$.
Then $\m_2(e_p,e_q)=e^{-\alpha_1}\,[T'_1]$ is a nonzero multiple of the unit 
in $CF^*(T'_1,T'_1)$, and $\m_2(e_q,e_p)=e^{-\alpha_1}\,[T_2]$ 
is a nonzero multiple of the unit in $CF^*(T_2,T_2)$: this makes it
reasonable to state that $T'_1$ and $T_2$ are isomorphic.

This example illustrates the failure of convergent power series
Floer homology to be invariant under Hamiltonian isotopies, even without
wall-crossing (recall the isotopy from $T_1$ to $T'_1$ did not cross
$f^{-1}(0)$); this is of course very different from the situation over
the Novikov ring. When we deform $T'_1$ back to $T_1$, we end up being
able to cancel all the intersection points even though they represent
nontrivial elements in Floer homology, because the cancellations in the
Floer differential occur between families of discs with different symplectic
areas (something which wouldn't be possible over Novikov coefficients). 
At the critical instant in the deformation, the discs with area
$\alpha_1$ have shrunk to points, while the discs with area $\alpha_2$
become pinched annuli. At the end of the deformation, the tori $T_1$ and 
$T_2$ are disjoint, and the discs have become holomorphic annuli with
boundary in $T_1\cup T_2$. 

It would be tempting to hope that a souped up version of Floer theory
that also includes holomorphic annuli would be better behaved. However,
in that case we would immediately hit a divergence issue when working
with complex coefficients: indeed, there are $2^k$ families of holomorphic
annuli with boundary in $T_1\cup T_2$ which cover $k$-to-1 the annulus 
bounded by the circles $\gamma_1$ and $\gamma_2$ in $\C$.

Even without considering annuli, divergence issues are already responsible 
for the bad properties of convergent power series Floer homology exhibited
here -- first and foremost, the lack of invariance under the Hamiltonian
isotopy from $T_1$ to $T'_1$. Denote by $H:[0,1]\times X^0\to \R$ a
family of Hamiltonians whose time 1 flow sends $T_1$ to $T'_1$, and recall
that continuation maps $\Phi:CF^*(T_1,T_2)\to CF^*(T'_1,T_2)$ and 
$\bar\Phi:CF^*(T'_1,T_2)\to CF^*(T_1,T_2)$ can
be defined by counting index 0 solutions of perturbed holomorphic 
curve equations of the form 
\begin{equation}\label{eq:perthol}
\frac{\partial u}{\partial s}+J\left(\frac{\partial u}{\partial t}-
\chi(s)X_H(t,u(s,t))\right)=0.\end{equation}
Here $u:\R\times [0,1]\to X^0$ is a map with $u(\R\times \{0\})\subset T_1$
and $u(\R\times\{1\})\subset T_2$ and satisfying suitable asymptotic
conditions at infinity, $X_H$ is the Hamiltonian vector field associated to
$H$, and $\chi:\R\to [0,1]$ is a suitable smooth cut-off function.

In our case, $\Phi$ and $\bar\Phi$ are obviously zero since
$CF^*(T_1,T_2)=0$; this of course prevents $\Phi\circ\bar\Phi:
CF^*(T'_1,T_2)\to CF^*(T'_1,T_2)$ from being homotopic to identity
as expected. Specifically, the homotopy would normally be constructed
by considering exceptional index $-1$ solutions to (\ref{eq:perthol})
where the cut-off $\chi$ is equal to 1 near $\pm \infty$ and $\int_{\R}
(1-\chi)$ varies between $0$ and infinity. In the present case, a
calculation shows that
that there are infinitely many exceptional solutions -- in fact there
are $2^k$ solutions of energy $k(\alpha_2-\alpha_1)$ for each integer $k$,
which makes the homotopy divergent. (To see this,
choose the Hamiltonian isotopy from $T_1$ to $T'_1$ to be lifted from the
complex plane by the projection $f$, and look at similar continuation
maps between $CF^*(\gamma_1,\gamma_2)=0$ and $CF^*(\gamma'_1,\gamma_2)$
inside $\C\setminus \{\epsilon\}$. In that case, an explicit calculation
shows that there is an infinite
sequence of exceptional index $-1$ solutions to (\ref{eq:perthol}), 
wrapping once, twice, etc.\ around the annulus bounded by $\gamma_1$ and
$\gamma_2$. Moreover, the exceptional trajectory which wraps $k$ times around the
annulus in $\C\setminus\{\epsilon\}$ can be shown to admit $2^k$
$S^1$-families of lifts to $X^0$.)

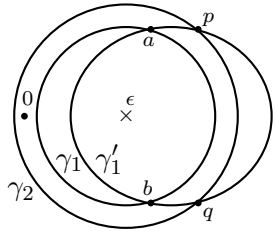
\begin{figure}[b]
\setlength{\unitlength}{1.5cm}
\begin{picture}(2,2)(-1,-1)
\psset{unit=\unitlength}
\pscircle(0,0){1}
\pscircle(0,0){0.8}
\psellipse(0.4,0)(0.9,0.8)
\pscircle*(-0.9,0){0.03}
\put(0,0){\makebox(0,0)[cc]{\tiny $\times$}}
\put(-0.92,0.1){\tiny $0$}
\put(0,0.1){\tiny $\epsilon$}
\put(-0.27,-0.45){\small $\gamma'_1$}
\put(-0.62,-0.45){\small $\gamma_1$}
\put(-1.05,-0.7){\small $\gamma_2$}
\pscircle*(0.64,-0.77){0.03} \put(0.68,-0.9){\tiny $q$}
\pscircle*(0.64,0.77){0.03} \put(0.68,0.83){\tiny $p$}
\pscircle*(0.22,-0.77){0.03} \put(0.15,-0.7){\tiny $b$}
\pscircle*(0.22,0.77){0.03} \put(0.15,0.62){\tiny $a$}
\end{picture}

\caption{$T_1$, $T'_1$ and $T_2$} \label{fig:wiggle2}
\end{figure}

Another instance of divergence occurs if we try to test the associativity
of the product in Floer homology. Namely, in addition to the isomorphisms
$e_p\in CF^0(T'_1,T_2)$ and $e_q\in CF^0(T_2,T'_1)$ considered above,
denote by $e_a\in CF^0(T'_1,T_1)$, resp.\ $e_b\in CF^0(T_1,T'_1)$, the
generators which come from the intersections in $f^{-1}(a)$, resp.\
$f^{-1}(b)$ (see Figure \ref{fig:wiggle2}). One easily checks that
$\m_2(e_a,e_b)$ is a nonzero multiple of the unit in $CF^*(T'_1,T'_1)$.
Then we can try to compose $e_a$, $e_b$ and $e_p$ in two different ways:
$\m_2(\m_2(e_a,e_b),e_p)$ is a nonzero multiple of $e_p$, whereas
$\m_2(e_a,\m_2(e_b,e_p))$ is zero since $\m_2(e_b,e_p)\in CF^*(T_1,T_2)=0$.
Passing to cohomology classes, this contradicts the expected associativity
of the product on Floer homology. A closer inspection reveals that this
is caused by the divergence of quantities such as $\m_3(e_a,f_a,e_p)$
(where $f_a$ is the generator of $CF^1(T_1,T'_1)$ corresponding to the
intersections in $f^{-1}(a)$): indeed, this triple product counts discs
obtained by cutting open the divergent series of annuli with boundary 
in $T_1\cup T_2$ already mentioned above.

In conclusion, there are many pitfalls associated to the use of convergent
power series Floer homology, even in fairly simple situations (compactifying
the above example to $\CP^2$, we would still encounter the same divergence
phenomena in a smooth projective Fano variety). A cautious view of the situation would dictate
that outside of the very simplest cases it is illusory to even attempt 
to work over complex coefficients, and that in general mirror symmetry 
is only a perturbative phenomenon taking place over a formal neighborhood 
of the large volume limit. Nonetheless, as long as one restricts oneself
to consider only certain aspects of Floer theory, the power series obtained
by working over the Novikov ring seem to often have good enough convergence
properties to allow the construction of a mirror that is an honest 
complex manifold (rather than a scheme over the Novikov field). Floer theory for
a single weakly unobstructed Lagrangian seems to be less prone to
divergence than the theory for pairs such as $(L_1,L_2)$ in the above example.
Also, in the example we have considered, divergence issues can be avoided
by equipping all our Lagrangian submanifolds with suitable Hamiltonian 
perturbation data (i.e., ``wiggling'' Lagrangians so that they intersect
sufficiently). However, more sophisticated divergent examples can be 
built e.g.\ inside conic bundles over elliptic curves; in some of these
examples, Floer products are given by series in $\Lambda_0$ for which
the radius of convergence is strictly less than 1, i.e.\ convergence
only holds for sufficiently large symplectic forms, regardless of
Hamiltonian perturbations.

\section{Relative mirror symmetry} \label{s:relms}

\subsection{Mirror symmetry for pairs}

In this section, we turn to mirror symmetry for a {\it pair} $(X,D)$, where
$X$ is a K\"ahler manifold and $D$ is a {\it smooth} Calabi-Yau hypersurface
in the anticanonical linear system. Our goal is to clarify the folklore
statement that ``the fiber of the mirror superpotential $W:X^\vee\to\C$ is
mirror to $D$''. The discussion is fairly similar to that in \S 7 of 
\cite{Au}.

Let $D\subset X$ be a hypersurface in the anticanonical linear system,
with defining section $\sigma\in H^0(X,K_X^{-1})$: then the
holomorphic volume form $\Omega=\sigma^{-1}\in \Omega^{n,0}(X\setminus D)$
(with poles along $D$) induces a holomorphic volume form $\Omega_D$ on $D$,
the {\it residue}
of $\Omega$ along $D$, characterized by the property that
$\Omega=\Omega_D\wedge d\log\sigma+O(1)$ in a neighborhood of $D$.
Additionally, the K\"ahler
form $\omega$ induces a K\"ahler form on $D$ by restriction.

It is reasonable to expect that special Lagrangian torus fibrations
on $X\setminus D$ should have a ``nice'' boundary behavior. Namely,
assuming that the K\"ahler metric on $X$ is complete, 
for a reasonable special Lagrangian fibration $\pi:X\setminus D\to B$
we expect:

\begin{conj}\label{conj:boundary}
Near $\partial B$, the fibers of $\pi$ are contained in a
neighborhood of $D$, and the smooth fibers are $S^1$-bundles over special Lagrangian tori
in $(D,\omega_{|D},\Omega_D)$.
\end{conj}

\noindent
(Here, by $\partial B$ we mean the part of the boundary 
of $B$ which lies at finite distance in the symplectic affine structure).

In other terms, we expect that near $D$ the special Lagrangian tori in
$X\setminus D$ accumulate onto special Lagrangian tori in $D$ (as observed
in the various examples we have discussed).
If Conjecture \ref{conj:boundary} holds, then $\partial B$ is the base
of a special Lagrangian fibration on $D$, and the (uncorrected) SYZ
mirror to $D$, $M_D$, can be identified as a complex hypersurface lying
inside the boundary of the (uncorrected) moduli space $M$ of pairs
$(L,\nabla)$ in $X\setminus D$.

Assume $D$ is smooth, and consider a special Lagrangian torus fiber
$L=\pi^{-1}(b)$ near $\partial B$: then we expect that $L$ bounds a 
distinguished family of Maslov index 2 holomorphic discs, namely small 
{\it meridian} discs in the normal direction to $D$. More precisely,
as $b$ approaches the boundary of $B$, we expect $L$ to collapse onto a special
Lagrangian torus $\Lambda$ in $D$, and the meridian discs to be approximated
by small discs inside the fibers of the normal bundle of $D$ lying above the points
of $\Lambda$.

Call $\delta$ the relative homotopy class of the meridian discs, and by
$z_\delta$ the corresponding holomorphic coordinate on $M$ (which is
also the contribution of the family of meridian discs to the superpotential).
Then we expect that $z_\delta$ is the dominant term in the superpotential
near the boundary of $M$, as the meridian discs have areas
tending to zero and all the other holomorphic discs have comparatively
much greater areas.

The boundary of $M$ corresponds to limiting pairs $(L,\nabla)$ where
the area of the meridian disc becomes $0$ (i.e., $L$ is entirely collapsed
onto a special Lagrangian torus in $D$); recalling that
$|z_\delta|=\exp(-\int_\delta\omega)$, this corresponds to $|z_\delta|=1$.
In fact, the boundary of $M$ fibers above the unit circle, via the map
\begin{equation}\label{eq:fibd}
\arg(z_\delta):\partial M=\{|z_\delta|=1\}\to S^1,
\end{equation}
with fiber $M_D=\{z_\delta=1\}$. The points of $M_D$ correspond to pairs
$(L,\nabla)$ where $L$ is entirely collapsed onto a special Lagrangian
torus $\Lambda\subset D$, {\it and} the holonomy of $\nabla$ around the
meridian loop $\mu=\partial\delta$ is trivial, i.e.\ $\nabla$ is pulled back from a $U(1)$ local
system on $\Lambda$. Thus $M_D$ is precisely the uncorrected SYZ mirror to
$D$.

In general, the fibration (\ref{eq:fibd}) has monodromy. Indeed, a local
trivialization is given by fixing a {\it framing}, i.e.\ an
$(n-1)$-dimensional subspace of $H_1(L,\Z)$ which under the projection
$L\to\Lambda$ maps isomorphically onto $H_1(\Lambda,\Z)$. (Less
intrinsically, we can choose a set of {\it longitudes},
i.e.\ lifts to $L$ of a collection of $n-1$ loops generating
$H_1(\Lambda,\Z)$); the framing data allows us to lift to $M$ a set of local
holomorphic coordinates on $M_D$. However, unless the normal bundle to $D$
is trivial there is no consistent {\it global} choice of framings: if we
move $\Lambda$ around a loop in $\partial B$ and keep track of a longitude
$\lambda$ lifting a loop $\gamma\in \Lambda$, the monodromy action is
of the form $\lambda\mapsto \lambda+k_\gamma\mu$, where $k_\gamma$ is the degree of the
normal bundle of $D$ over the surface traced out by $\gamma$. 

A more thorough calculation shows that the monodromy of (\ref{eq:fibd}) 
is given by a symplectomorphism of $M_D$ which geometrically realizes
(as a fiberwise translation in the special Lagrangian fibration $M_D\to
\partial B$ dual to the SYZ fibration on $D$) the mirror to the
autoequivalence $-\otimes K_{X|D}^{-1}$ of $D^b\mathrm{Coh}(D)$.

This is easiest to see if we
assume that, in a neighborhood of $D$, the anticanonical bundle $K_X^{-1}$
can be equipped with a semi-flat connection, i.e.\ a holomorphic connection
whose restriction to the fibers of $\pi$ is flat.
Then the parallel transport from one fiber of (\ref{eq:fibd}) to another
can be realized geometrically as follows:
given a pair $(L,\nabla)$ where $L$ is almost
collapsed onto a special Lagrangian $\Lambda\subset D$, we can modify the
holonomy of $\nabla$ around the meridian loop by adding to it a multiple
of $\mathrm{Im}(\sigma^{-1}\partial\sigma)_{|L}$, where $\sigma$ is the defining section of
$D$. The monodromy is then $(L,\nabla)\mapsto
(L,\nabla+\mathrm{Im}(\sigma^{-1}\partial\sigma)_{|L})$, which in the limit where
$L$ collapses onto $D$ is exactly the expected transformation.%
\medskip

If we can neglect the terms other than $z_\delta$ in the superpotential,
for instance in the {\it large volume limit}, then $M_D$ is essentially
identified with the fiber of $W$ at 1. In fact, recall from the
discussion at the end of \S \ref{ss:lg} that
changing the K\"ahler class to $[\omega]+tc_1(X)$ ``enlarges'' the mirror
while rescaling the superpotential by a factor of $e^{-t}$: thus, assuming
that $X$ is Fano, or more generally that $-K_X$ is nef, the flow to the
large volume limit can be realized simply by rescaling the superpotential.
Hence, Conjecture \ref{conj:boundary} implies:

\begin{conj}\label{conj:fiber}
If $(X^\vee,W)$ is mirror to $X$, and if $-K_X$ is nef, then
for $t\to \infty$
the family of hypersurfaces $\{W=e^t\}\subset X^\vee$ is asymptotic (up to
corrections that decrease exponentially with $t$) to the
family of mirrors to $(D,\omega_{|D}+tc_1(X)_{|D})$.
\end{conj}

For example, considering the mirror to $\CP^2$ with
$[\omega]\cdot[\CP^1]=\Lambda$, the $j$-invariant of the elliptic curve
$\{x+y+e^{-\Lambda}/xy=e^t\}\subset (\C^*)^2$ can be determined to equal
$$\frac{e^{3t+\Lambda}(e^{3t+\Lambda}-24)^3}{e^{3t+\Lambda}-27}=
e^{9t+3\Lambda}+\dots,$$ whose leading term matches with the symplectic
area of the anticanonical divisor after inflation (observe that
$([\omega]+tc_1)\cdot [\CP^1]=3t+\Lambda$).

There are two reasons why this statement only holds asymptotically for
$t\to\infty$. First, the formula for the superpotential includes other
terms besides $z_\delta$, so the hypersurfaces $\{W=e^t\}$ and
$\{z_\delta=e^t\}$ are not quite the same. More importantly, the instanton
corrections to the mirror of $D$ are {\it not} the same as the instanton
corrections to the fiber of $z_\delta$. When constructing the
mirror to $X$, the geometry of $M_D\subset M$ gets corrected by 
wall-crossing terms that record holomorphic Maslov index 0 discs in $X$; 
whereas, when constructing the mirror of $D$, the corrections only arise
from Maslov index 0 holomorphic discs in $D$. 

In other terms: the instanton corrections to the mirror of $X$ arise from
walls generated by singularities in the fibration $\pi:X\setminus D\to B$
(i.e., singularities in the affine structure of $B$),
whereas the instanton corrections to the mirror of $D$ arise from the
walls generated by singularities in the fibration $\pi_D:D\to \partial B$
(i.e., singularities in the affine structure of $\partial B$).
Since the singularities of the affine structure on $\partial B$ are 
induced by those strata of singularities of $B$ that hit the boundary,
the wall-crossing phenomena in $D$ are induced
by a subset of the wall-crossing phenomena 
in $X$, but there are also walls in $X$ which hit the boundary of $B$ without
being induced by singularities at the boundary.

On the other hand, the smooth fibers of $W$ are symplectomorphic to each
other and to the hypersurface $\{z_\delta=1\}$. Moreover, it is
generally believed that the K\"ahler class of the mirror should not be 
affected by instanton corrections, so the discrepancy discussed above is no
longer an issue. Hence: we expect that the fibers of $W$, viewed 
{\it as symplectic
manifolds}, are mirror to the divisor $D$ viewed {\it as a complex manifold}.
(Observe that, from this perspective, the parameter $t$ in Conjecture
\ref{conj:fiber} no longer plays any role, and accordingly the geometries
are expected to match on the nose.)


\subsection{Homological mirror symmetry}\label{ss:relhms}

Assuming Conjectures \ref{conj:boundary} and \ref{conj:fiber}, we can
try to compare the statements of homological mirror symmetry for $X$
and for the Calabi-Yau hypersurface $D$. Due to the mismatch between
the complex structure on the mirror to $D$ and that on the fibers of $W$
(see Conjecture \ref{conj:fiber}), in general we can only hope to achieve this
in one direction, namely relating the derived categories of coherent
sheaves on $X$ and $D$ with the Fukaya categories of their mirrors.

Denote by $(X^\vee,W)$ the
mirror to $X$, and by $D^\vee$ the mirror to $D$, which we identify
symplectically with a fiber of $W$, say $D^\vee=\{W=e^t\}\subset X^\vee$
for fixed $t\gg 0$. 

First we need to briefly describe the Fukaya category
of the Landau-Ginzburg model $W:X^\vee\to\C$. The general idea, which
goes back to Kontsevich \cite{KoENS} and Hori-Iqbal-Vafa \cite{HIV}, is
to allow as objects {\em admissible Lagrangian submanifolds}\/ of $X^\vee$;
these can be described either as
potentially non-compact Lagrangian submanifolds which, outside of a compact
subset, are invariant under the gradient flow of $-\mathrm{Re}(W)$, or,
truncating, as compact Lagrangian submanifolds
with (possibly empty) boundary contained inside a fixed reference fiber
of $W$ (and satisfying an additional condition). The case
of Lefschetz fibrations (i.e., when the critical points of $W$ are
nondegenerate) has been studied in great detail by Seidel; in this case,
which is by far the best understood, the theory can be
formulated in terms of the vanishing cycles at the critical points 
(see e.g.\ \cite{SeBook}).

The formulation which is the most relevant to us is
the one discussed by Abouzaid in \cite{abouzaid}: in this version,
one considers Lagrangian submanifolds of $X^\vee$ with boundary contained 
in the given reference fiber $D^\vee=W^{-1}(e^t)$,
and which near the reference fiber 
are mapped by $W$ to an embedded curve $\gamma\subset\C$. 

\begin{definition}
A Lagrangian submanifold $L\subset X^\vee$ with (possibly empty)
boundary $\partial L\subset D^\vee=W^{-1}(e^t)$
is {\em admissible with phase $\varphi\in(-\frac\pi2,\frac\pi2)$} 
if\/ $|W|<e^t$ at every point of $\mathrm{int}(L)$ and, 
near $\partial L$, the restriction of $W$ to $L$
takes values in the half-line $e^t-e^{i\varphi}\R_+$.
\end{definition}

Floer theory is then defined by choosing a specific set of Hamiltonian
perturbations, which amounts to deforming the given admissible
Lagrangians so that their phases are in increasing order, and ignoring
boundary intersections. For instance, to determine $HF(L_1,L_2)$, one first
deforms $L_2$ (rel.\ its boundary) to an admissible Lagrangian $L_2^+$ whose
phase is
greater than that of $L_1$, and one computes Floer homology for the pair of
Lagrangians $(L_1,L_2^+)$ inside $X^\vee$, ignoring boundary intersections.
We denote by $\mathcal{F}(X^\vee,D^\vee)$ the Fukaya category constructed
in this manner. (In fact, strictly speaking, one should place the reference
fiber ``at infinity'', i.e.\ either consider a
limit of this construction as $t\to +\infty$, or enlarge the symplectic
structure on the subset $\{|W|<e^t\}$ of $X^\vee$ so that the symplectic
form blows up near the boundary and the K\"ahler metric becomes complete;
for simplicity we ignore this subtlety.)


By construction, the boundary of an admissible Lagrangian in $X^\vee$
is a Lagrangian submanifold
of $D^\vee$ (possibly empty, and not necessarily connected).
There is a {\it restriction $A_\infty$-functor}
$\rho:\mathcal{F}(X^\vee,D^\vee)\to \mathcal{F}(D^\vee)$ from the 
Fukaya category of the Landau-Ginzburg model $(X^\vee,W)$ to the
(usual) Fukaya category of $D^\vee$. At the level of objects,
this is simply $(L,\nabla)\mapsto (\partial
L,\nabla_{|\partial L})$. At the level of morphisms, the $A_\infty$-functor
$\rho$ consists of a collection
of maps $$\rho_{(k)}:\mathrm{Hom}_{\mathcal{F}(X^\vee,D^\vee)}(L_1,L_2)\otimes
\dots\otimes \mathrm{Hom}_{\mathcal{F}(X^\vee,D^\vee)}(L_k,L_{k+1})\to
\mathrm{Hom}_{\mathcal{F}(D^\vee)}(\partial L_1,\partial L_{k+1}).$$
The first order term $\rho_{(1)}$ is the easiest to describe: given an intersection
point $p\in \mathrm{int}(L_1)\cap \mathrm{int}(L_2^+)$, $\rho_{(1)}(p)$ is a
linear combination of intersection points in which the coefficient of
$q\in \partial L_1\cap \partial L_2$ counts holomorphic strips
in $(X^\vee,L_1\cup L_2^+)$ connecting $p$ to $q$. Similarly, given $k+1$
admissible Lagrangians $L_1,\dots,L_{k+1}$, and perturbing them so that 
their phases are in increasing order, $\rho_{(k)}$ counts holomorphic
discs in $(X^\vee,\bigcup L_i^+)$ with $k$ corners at prescribed interior intersection points and
one corner at a boundary intersection point.

Homological mirror symmetry for the pair $(X,D)$ can then be summarized
by the following conjecture:

\begin{conj}\label{conj:relhms}
There is a commutative diagram
$$\begin{CD}
D^bCoh(X) @>\text{restr}>> D^bCoh(D) \\
@V{\simeq}VV @VV{\simeq}V \\
D^\pi\mathcal{F}(X^\vee,D^\vee) @>{\rho}>> D^\pi\mathcal{F}(D^\vee)
\end{CD}
$$
\end{conj}

In this diagram, the horizontal arrows are the restriction functors, and
the vertical arrows are the equivalences predicted by homological mirror
symmetry. The reader is referred to \cite{AKO2} for a verification in the
case of Del Pezzo surfaces.
\medskip

Another type of Fukaya category that can be associated to $X^\vee$ is
its {\it wrapped Fukaya category} $\mathcal{F}_{wr}(X^\vee)$ \cite{AS}. The objects of that category
are again non-compact Lagrangian submanifolds, but the Hamiltonian 
perturbations used to define Floer homology now diverge at infinity.
Assuming that $W$ is proper, we can e.g.\ use the Hamiltonian flow generated
by a function of $|W|$ that increases sufficiently quickly at infinity;
however, the wrapped category can be defined purely in terms of the
symplectic geometry of $X^\vee$ at infinity, without reference to the
superpotential (see \cite{AS}).

Homological mirror symmetry for the open Calabi-Yau $X\setminus D$ then
predicts an equivalence between the derived category of coherent sheaves
$D^bCoh(X\setminus D)$ and the derived wrapped Fukaya category 
$D^\pi\mathcal{F}_{wr}(X^\vee)$. Moreover, the restriction functor from
$D^bCoh(X)$ to $D^bCoh(X\setminus D)$ is expected to correspond to a
natural functor $\varpi$ from the Fukaya category of the Landau-Ginzburg model
$(X^\vee,W)$ to the wrapped Fukaya category of $X^\vee$. On objects, 
$\varpi$ is essentially identity (after sending the reference fiber to infinity,
or extending admissible Lagrangians to non-compact ones by parallel transport
along the gradient flow of $\mathrm{Re}(W)$). On morphisms, $\varpi$ is
essentially an inclusion map if we set up the Hamiltonian perturbations
in the wrapped category to be supported outside of the region where
$|W|<e^t$; or, more intrinsically, $\varpi$ is the continuation
map induced on Floer complexes by the deformation from the small
Hamiltonian perturbations used to define the Fukaya category of $(X^\vee,W)$
to the large Hamiltonian perturbations used to define the wrapped category.

In fact, the wrapped Fukaya category can alternatively be defined from
$\mathcal{F}(X^\vee,W)$ as the result of localization with respect to a
certain natural transformation from the Serre functor (up to a shift) to
the identity, induced by the monodromy of $W$ near infinity (see \S 4 of
\cite{SeHoch} and \S 6 of \cite{SeHoch2}); this parallels
the fact that $D^bCoh(X\setminus D)$ is the localization of $D^bCoh(X)$ with respect
to the natural transformation from $-\otimes K_X$ (i.e., the Serre 
functor up to a shift) to the identity given by the defining section of~$D$.

Finally, when considering compact closed Lagrangian submanifolds, there
is no difference between the Fukaya category of $(X^\vee,W)$ and the wrapped
Fulaya category; the full subcategory consisting of these compact objects
is expected to be equivalent to the subcategory of $D^bCoh(X\setminus D)$
generated by complexes with compactly supported cohomology.

\subsection{Complete intersections}

As pointed out to the author by Ludmil Katzarkov, the above ideas can
be extended to understand mirror symmetry for complete
intersections (remaining in the framework of manifolds with effective
anticanonical divisors). Namely, consider 
divisors $D_1,\dots,D_k\subset X$ (smooth, or at most with normal crossing
singularities), intersecting each other
transversely, such that $\sum D_i=-K_X$. Let $(X^\vee,W)$ be the
mirror of $X$ relative to the anticanonical divisor $\sum D_i$:
then the superpotential on $X^\vee$ splits into a sum $W=W_1+\dots+W_k$,
where $W_i:X^\vee\to\C$ records the contributions to the superpotential
of holomorphic Maslov index 2 discs which hit the component $D_i$ of the
anticanonical divisor.

For a subset $I\subseteq \{1,\dots,k\}$, consider the complete intersection
$X_I=\bigcap_{i\in I} D_i\subset X$, and the divisors $D_{I,j}=X_I\cap D_j$,
$j\not\in I$, whose sum represents the
anticanonical class of $X_I$. Then we have:

\begin{conj}\label{conj:ci}
In the large volume limit $t\to \infty$, the mirror to $X_I$ equipped
with the K\"ahler form $\omega_{|X_I}+tc_1(X)_{|X_I}$ and
the anticanonical divisor $\sum_{j\not\in I} D_{I,j}$ is approximated
(in the sense of Conjecture \ref{conj:fiber}) by
the complete intersection 
$X^\vee_I:=\bigcap_{i\in I} W_i^{-1}(e^t)$ in $X^\vee$, equipped with
the superpotential $W_I:=\sum_{j\not\in I} W_j$.
\end{conj}

As before, if we are only interested in comparing the complex geometry
of $X_I$ with the symplectic geometry of $(X^\vee_I,W_I)$, then the
construction does not depend on the parameter $t$, and passage
to the large volume limit is not needed.

Conjecture \ref{conj:ci} can be understood geometrically as follows.
In this setting, we expect to have a special Lagrangian torus fibration
$\pi:X\setminus (\bigcup D_i)\to B$, whose base $B$ has boundary and
corners: at the boundary, the
special Lagrangian fibers collapse onto one of the hypersurfaces $D_i$,
and at the corners they collapse onto the
intersection of several $D_i$. (This picture is e.g.\ obvious in the
toric setting, where $B$ is the interior of the moment polytope.)

Whenever the fibers of $\pi$ lie sufficiently close to $D_i$, they are
expected to bound small meridian discs intersecting $D_i$ transversely
once, whereas the other families of discs have comparatively larger
symplectic area, so that $W_i=z_{\delta_i}+o(1)$. Setting $z_{\delta_i}$ equal to 1
for $i\in I$ amounts to considering special Lagrangian tori that are
completely collapsed onto $X_I=\cap_{i\in I}D_i$, equipped with flat
connections that have trivial holonomy along the meridian loops, i.e.\
are pulled back from special Lagrangian tori in $X_I$. Thus, before
instanton corrections, $\bigcap_{i\in I} \{z_{\delta_i}=1\}$ is the
(uncorrected) SYZ mirror to $X_I\setminus (\bigcup_{j\not\in I} D_{I,j})$.
When $t\to \infty$ the discrepancy between $W_i$ and $z_{\delta_i}$ and
the differences in instanton corrections are expected to become negligible.

Moreover, in the limit where $L\subset X\setminus (\bigcup D_i)$ 
collapses onto a special Lagrangian $\Lambda\subset X_I\setminus
(\bigcup_{j\not\in I} D_{I,j})$, for $j\not\in I$ the dominant terms 
in $W_j$ should correspond to families of holomorphic discs in $(X,L)$ 
that converge to holomorphic discs in $(X_I,\Lambda)$ (intersecting $D_{I,j}$).
Hence, $\sum_{j\not\in I} W_j$ should differ from the superpotential for
the mirror to $X_I$ by terms that become negligible in the large volume
limit.

As a special case of Conjecture \ref{conj:ci}, taking $I=\{1,\dots,k\}$,
(in the large volume
limit) the fiber of $(W_1,\dots,W_k)$ is mirror to the 
Calabi-Yau complete intersection
$X_{\{1,\dots,k\}}=D_1\cap\dots\cap D_k$. (In this case there is no 
residual superpotential.) This is consistent with standard conjectures.
\medskip

It is also worth noting that, in a degenerate toric limit, Conjecture
\ref{conj:ci} recovers the predictions made by Hori and Vafa \cite{HV}
for mirrors of Fano complete intersections in toric varieties.
To give a simple example, consider $X=\CP^3$ (with
$\int_{\CP^1}\omega=\Lambda$), and let $D_1,D_2\subset X$ be quadric
surfaces intersecting transversely in an elliptic curve $E=D_1\cap D_2$.
Then the superpotential on $X^\vee$ decomposes as a sum $W=W_1+W_2$.
In the degenerate limit where $D_1$ and $D_2$ are {\it toric} quadrics
consisting of two coordinate hyperplanes each, and $E$ is a singular
elliptic curve with four rational components, we have
$X^\vee=\{z_0z_1z_2z_3=e^{-\Lambda}\}\subset (\C^*)^4$, and 
$W=W_1+W_2$, where $W_1=z_0+z_1$ and $W_2=z_2+z_3$.
Then the mirror to $D_1$ is the surface $$\{z_0z_1z_2z_3=e^{-\Lambda},
\ z_0+z_1=e^t\}\subset (\C^*)^4,$$ equipped with the superpotential $W_2=z_2+z_3$,
and similarly for $D_2$; and the mirror to $E$ is the curve $\{z_0z_1z_2z_3=
e^{-\Lambda},\ z_0+z_1=e^t,\ z_2+z_3=e^t\}$ 
(a noncompact elliptic curve with four punctures). These
formulas are essentially identical to those in Hori-Vafa \cite{HV}.
To be more precise: viewing $D_i$ and $E$ as symplectic manifolds (in
which case the degeneration to the toric setting should be essentially
irrelevant, i.e.\ up to a fiberwise compactification of the Landau-Ginzburg
models we can think of smooth quadrics and elliptic curves), but taking
the {\it large volume limit}\/ $t\to \infty$,
these formulas give an {\it approximation} to the complex geometry of
the mirrors. On the other hand, if we consider the symplectic geometry of
the mirrors, then the formulas give exact mirrors to $D_i$ and $E$ viewed
as singular complex manifolds
(torically degenerated quadrics and elliptic curves, i.e.\ large complex structure limits).
Thus Hori and Vafa's formulas for toric complete intersections
should be understood as a construction
of the mirror at a limit point in {\it both} the complex and K\"ahler moduli
spaces.

\end{document}